\documentclass{article}

\usepackage{microtype}
\usepackage{graphicx}
\usepackage{subfigure}
\usepackage{multirow}
\usepackage{float}
\usepackage{soul}
\usepackage{booktabs} 
\usepackage{hyperref}
\usepackage{longtable}
\usepackage{lscape}
\usepackage{caption}

\usepackage[accepted]{icml2024}

\usepackage{amsmath}
\usepackage{amssymb}
\usepackage{mathtools}
\usepackage{amsthm}
\usepackage{amsfonts,mathrsfs}
\usepackage{enumitem}
\usepackage[capitalize,noabbrev]{cleveref}

\theoremstyle{plain}

\theoremstyle{definition}

\theoremstyle{remark}

\usepackage[textsize=tiny]{todonotes}

\begin{document}

\twocolumn[

\icmltitle{LLaMoCo: Instruction Tuning of Large Language Models for \\Optimization Code Generation}


\icmlsetsymbol{corresponding}{*}

\begin{icmlauthorlist}
\icmlauthor{Zeyuan Ma}{scut}
\icmlauthor{Hongshu Guo}{scut}
\icmlauthor{Jiacheng Chen}{scut}
\icmlauthor{Guojun Peng}{scut}\\
\icmlauthor{Zhiguang Cao}{smu}
\icmlauthor{Yining Ma}{ntu}
\icmlauthor{Yue-Jiao Gong}{scut}

\end{icmlauthorlist}

\icmlaffiliation{scut}{School of Computer Science and Engineering, South China University of Technology, Gungzhou, Guangdong, China}

\icmlaffiliation{smu}{School of Computing and Information Systems, Singapore Management University, Singapore.}

\icmlaffiliation{ntu}{Nanyang Technological University, Singapore}

\icmlcorrespondingauthor{Yining Ma}{yiningma@u.nus.edu}
\icmlcorrespondingauthor{Yue-Jiao Gong}{gongyuejiao@gmail.com}

\vskip 0.3in
]



\printAffiliationsAndNotice{}  

\begin{abstract}

Recent research explores optimization using large language models (LLMs) by either iteratively seeking next-step solutions from LLMs or directly prompting LLMs for an optimizer. However, these approaches exhibit inherent limitations, including low operational efficiency, high sensitivity to prompt design, and a lack of domain-specific knowledge. We introduce LLaMoCo, the first instruction-tuning framework designed to adapt LLMs for solving optimization problems in a code-to-code manner. Specifically, we establish a comprehensive instruction set containing well-described problem prompts and effective optimization codes. We then develop a novel two-phase learning strategy that incorporates a contrastive learning-based warm-up procedure before the instruction-tuning phase to enhance the convergence behavior during model fine-tuning. The experiment results demonstrate that a CodeGen (350M) model fine-tuned by our LLaMoCo achieves superior optimization performance compared to GPT-4 Turbo and the other competitors across both synthetic and realistic problem sets. The fine-tuned model and the usage instructions are available at \url{https://anonymous.4open.science/r/LLaMoCo-722A}.

\end{abstract}

\section{Introduction}
\label{sec:intro}
Nowadays, Large Language Models~(LLMs) are posing a profound impact on human society~\cite{floridi2020gpt,llm-app-2}. Through text generation, LLMs exhibit extraordinary prowess in natural language understanding and adeptness in solving complex tasks~\cite{llm-app-1,llm-app-2,llm-app-3}. This prompts a research question: Can LLMs even handle the challenging \textit{Optimization} problems that are usually difficult for humans to address? This forms the core of our study in this paper.


In the literature, several existing works have been developed to explore the possibilities of solving optimization problems using LLMs. A typical way is to iteratively prompt LLMs to output better solutions through a multi-turn conversation with LLMs~\cite{opro,opro-to,opro-ntu,opro-moo}. 
Typically, they leverage LLMs through an iterative process, sometimes incorporating the concept of in-context learning. This involves the steps of presenting the LLMs with a set of initial or current best-so-far solutions and iteratively requesting LLMs to generate potentially superior solutions.
While showing certain effectiveness in solving optimization tasks, these solution-to-solution approaches have several limitations: 1) the scale of target optimization tasks~(e.g., the number of variables, historical solutions and newly generated solutions) is \textit{constrained by the context window length of LLMs}; 2) the iterative process typically involves \textit{hundreds rounds of conversations}, consuming multitudinous resources; and 3)~due to the LLMs' \textit{sensitivity to prompt design}, it is nontrivial to provide coherent prompts for LLMs to guarantee ideal outputs.
\begin{figure*}[ht]
\vskip 0.1in
\centerline{\includegraphics[width=1.98\columnwidth]{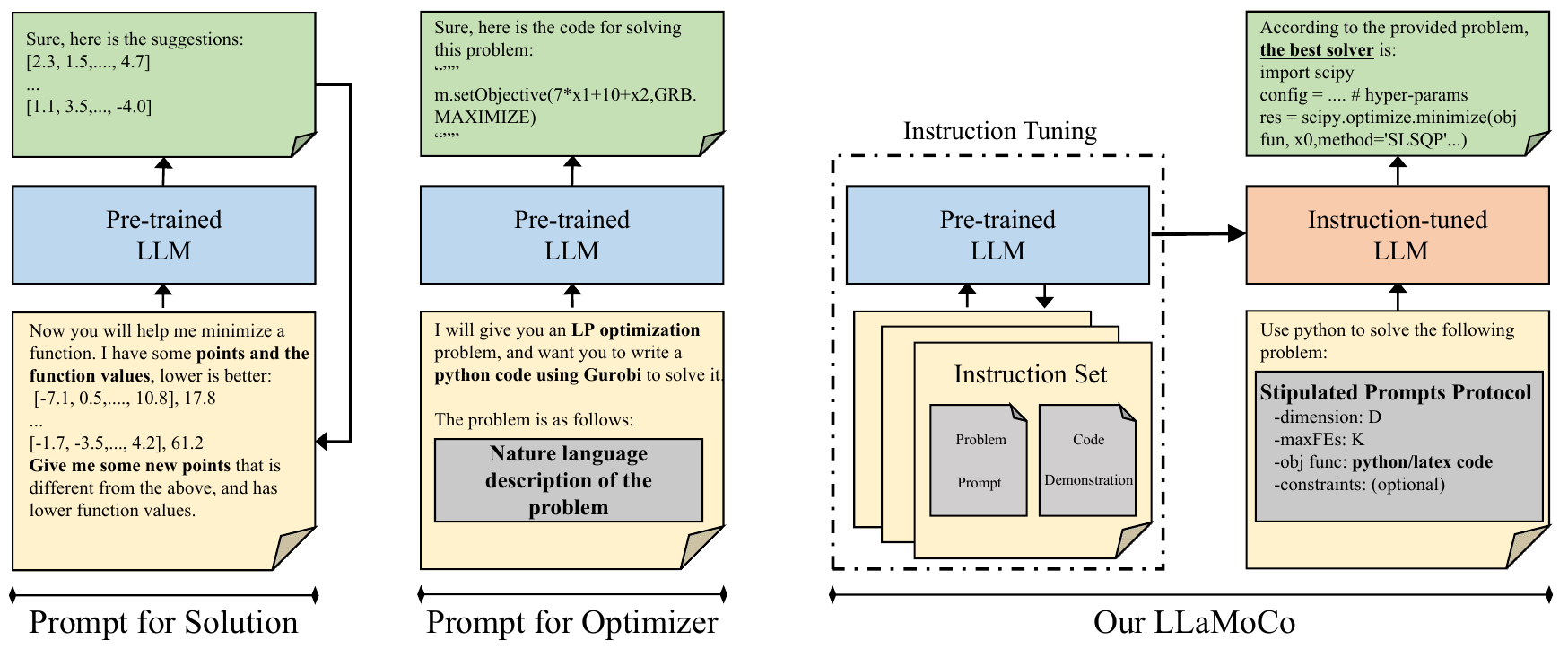}}

\caption{Conceptual overview of \textit{LLMs as optimizers}. \textbf{Left}: optimization through iteratively prompting LLMs for better solutions~(solution-to-solution style), such as OPRO~\cite{opro}. 
\textbf{Middle}: optimization through directly prompting LLMs for an optimizer with code implementation, such as OptiMus~\cite{optimus}. To ensure rational output, the prompts should include hints about the type of problem and the suggested optimizer.
\textbf{Right}: our LLaMoCo, which first tunes general LLMs on a problem-code instruction set, then can be used to generate proper optimization code given the formatted problem prompts.}
\label{overview1}
\end{figure*}

An alternative way is to directly prompt LLMs for optimization programs, namely a piece of executable code that either reuses existing optimization toolboxes~\cite{optimus} or combines multiple high-performing optimizers to create a novel one~\cite{opro-alg}. It could be more efficient than the solution-to-solution methods for two reasons: 1)~only a simple or few rounds of conversation are necessary for generating codes; and 2)~the prompts and the generated codes do not include solution information, making it compatible as the problem scales. However, it remains crucial to carefully craft prompts to ensure logical coherence of the generated codes. For example, OptiMus~\cite{optimus}  integrates hints about the optimizer to be generated directly into the prompts, necessitating a deep understanding of optimization techniques and expertise in the field. Additionally, 
using the LLMs pre-trained on a wide range of corpus currently falls short in generating a customized optimizer tailored to a specific optimization problem instance. This limitation is identified as the \textit{lack of domain-specific expert knowledge}~\cite{llm-survey}, which also extends to other intricate tasks with structured data, such as knowledge-base question answering and semantic parsing~\cite{structgpt,strcture}.

In this paper, we propose LLaMoCo, a novel framework that fine-tunes general-purpose \underline{L}arge \underline{La}nguage \underline{M}odels for \underline{o}ptimization \underline{Co}de generation. Different from the above approaches that are sorely based on prompt engineering, our LLaMoCo fine-tunes the LLMs on a well-formatted instruction set comprising code-to-code pairs of problem prompts and executable optimization programs. Once the training is completed, the fine-tuned model can be generalized to unseen optimization problems, i.e.,  crafting an optimizer based on the specific problem structure. Our LLaMoCo holds the following advantages against the preliminary works. 1) The solution information-free setting, where an optimization program is generated in a single round of conversation, makes it easier to handle large-scale problems with higher efficiency than the solution-to-solution methods. 2) The stipulated prompt protocol for users to describe their optimization problems minimizes the domain knowledge and efforts required for prompt design. 3) The fine-tuned LLMs by our LLaMoCo provide users with more robust and expert-level optimizers than those obtained by directly using general code-generating LLMs. In~\cref{overview1}, we illustrate the difference between LLaMoCo and existing approaches that leverage LLMs for optimization.

However, achieving expert-level LLMs for optimization tasks presents certain challenges. To overcome these challenges, we contribute to the following aspects: 1)~We establish the first instruction set for fine-tuning LLMs as expert-level optimizer generators.  This instruction set offers meticulously crafted problem descriptions and corresponding well-performing optimizer implementations selected from a wide spectrum of advanced optimizers, refined through extensive benchmarking with fine-grained hyper-parameter search; 2)~
We put forth a two-phase adaption strategy, which first enhances the latent space representation of a given problem instance through contrastive learning~\cite{cl}, followed by the conventional sequence-to-sequence loss for instruction tuning. Such design significantly accelerates the convergence of fine-tuned LLMs, resulting in superior performance; 3)~LLaMoCo has been meticulously designed for user-friendliness. Users can focus on the optimization problem itself following a stipulated prompt protocol, and then the prompt is automatically constructed and fed into the LLMs fine-tuned by our LLaMoCo.

Our benchmark experiments reveal the remarkably robust optimization performance of our LLaMoCo, surpassing existing methods. Notably, we show that instruction tuning of a relatively small LLM, e.g., CodeGen-$350$M~\cite{codegen}, on domain-specific tasks can yield substantial performance enhancements, even surpassing the very large and powerful models like GPT-$4$~\cite{gpt-4t}. Moreover, we provide in-depth analyses of the proposed two-phase adapting strategy, the sensitivity to training data distribution, and the zero-shot generalization performance. 

\textbf{In summary, our contributions are four folds}: 
1)~Introduction of LLaMoCo, the first instruction tuning framework for adapting general-purpose LLMs for generating expert-level optimizers. 2)~Establishment of the large-scale instruction set on optimization domain, providing copious code implementation of advanced optimizers at instance level~(\cref{sec3:instruction_set}). 3)~Development of a novel two-phase training strategy that reinforces the latent space representations of the prompts through efficient contrastive warm-up training, boosting the subsequent instruction tuning performance~(\cref{sec3:instruction_tuning}). 4)~Demonstration of LLaMoCo's superior optimization performance against existing LLM-based optimizers. The fine-tuned LLMs exhibit remarkable zero-shot generalization ability to realistic optimization tasks, with certain efficiency and code robustness~(\cref{sec4:exp}).





\section{Related Works}
\subsection{Fine-tuning LLMs}
Pre-trained Large Language Models (LLMs) can be refined by additional parameter updates on specific tasks through a fine-tuning process. We introduce two prominent fine-tuning strategies: Instruction Tuning~(IT)~\cite{IT1} and Alignment Tuning~(AT)~\cite{rlhf1,rlhf2}, each serving distinct purposes. Generally, IT involves fine-tuning pre-trained LLMs using a moderate collection of formatted task instances~\cite{IT2}. The fine-tuning process typically includes two steps: 1) prepare instruction-formatted training examples by associating a task description with each task instance, which aids LLMs in understanding tasks through the instructions~\cite{formatteddata}; and 2) leverage the prepared instruction set to fine-tune LLMs using a sequence-to-sequence supervised loss~\cite{gupta2023instruction}. By incorporating a well-established task-specific instruction set, IT can be an effective approach to inject domain-specific knowledge into general LLMs. This enables the transfer of LLMs to specific experts in domains like medicine~\cite{medicine-IT}, law~\cite{lawyer-IT} and finance~\cite{finance-IT}. 

Differently, AT aims to correct unexpected behaviors of LLMs by aligning the models with human values and preferences~\cite{IT1,rlhf2}. A practical algorithm for AT is the Reinforcement Learning from Human Feedback~(RLHF)~\cite{rlhf2}, which firstly estimates a reward model on a human-preference data collection via maximum likelihood. It then uses the learned reward model to provide feedback and post-trains the LLMs through Proximal Policy Optimization~(PPO)~\cite{ppo}. A recent work named Direct Preference Optimization~(DPO)~\cite{dpo} first reparameterizes the reward function based on the parameters of the pre-trained LLMs, saving the modelling and training of the reward function. DPO is mathematically equivalent to RLHF but is even more efficient, which is widely adopted in the latest LLMs such as Mistral 8x7B~\cite{mistral}.  

\subsection{LLMs for Code Generation}
The task of generating code from natural language descriptions is both exciting and inherently complex~\cite{llm-code-survey,chen2021evaluating}. 
Although general-purpose LLMs such as GPT~\cite{gpt}, Llama $2$~\cite{llama2} and Mistral~\cite{mistral} show competitive performance on the widely used LLM benchmarks including HumanEval~\cite{chen2021evaluating}, MBPP~\cite{mbpp} and DS-$1000$~\cite{ds-1000},  their performance on a particular task may still be limited. 
Recent efforts have focused on developing Large Language Models (LLMs) specifically tailored for code generation. These models are either trained exclusively on code, such as AlphaCode~\cite{alphacode} and StarCoder~\cite{starcoder}, or fine-tuned from general LLMs, like Codex~\cite{chen2021evaluating} and Code Llama~\cite{codellmma}. Notably, Codex shows that a $12$B LLM can solve $72.31\%$ of complex programming tasks posed by humans.  This success has led to the emergence of various Code LLMs, such as CodeGen~\cite{codegen} that factorizes a potentially long specification into multiple steps to enhance program synthesis, and Code Llama that increases Llama $2$ models through a cascade of fine-tuning steps. Other models such as Phi-$2$~\cite{phi-2}, InCoder~\cite{incoder} and CodeGeeX~\cite{codegeex} have also gained great attention.       
\subsection{LLMs as Optimizers}
Optimization plays a crucial role in numerous science and engineering fields but poses significant challenges. Unlike simpler tasks such as language understanding that can be easily handled by humans, optimization tasks can hardly be solved by humans without efficient algorithms. 
The underlying complexity of solving optimization problems tests the reasoning and generalization abilities of LLMs. Recently, there are several works that explore the potential use of LLMs as optimizers~\cite{opro,opro-alg,opro,guo2023connecting}, mostly based on prompt engineering and sometimes in an in-context learning way~\cite{icl}. Typically, these methods consider a set of candidate solutions to be improved. LLMs receive prompts containing these solutions and their objective values to propose improved ones. This process iterates until a termination condition is reached. Moreover, several studies introduce additional instructions, such as the mutation and crossover operations, to the naive prompts. This enables LLMs to mimic human-developed evolutionary operators, thereby achieving improved performance.~\cite{opro-ntu,opro-moo,lehman2023evolution,chen2023evoprompting}. However, these approaches have limitations in efficiency due to the need for extensive iterations. In contrast, several studies consider prompting LLMs directly for optimization programs, focusing on either creating new optimizers~\cite{opro-alg} or leveraging the combination of existing ones~\cite{optimus}. To the best of our knowledge, all the aforementioned works focus on prompt engineering of pre-trained LLMs, and the area of fine-tuning general LLMs with optimization-domain knowledge remains unexplored.

\section{LLaMoCo}
We introduce LLaMoCo, the first instruction tuning framework for adapting general-purpose LLMs as optimizers. It operates on a code-to-code basis, whereby, given anoptimization problem where its objective function and constraints are described using Python or LaTeX codes, the fine-tuned LLMs would generate a code implementation of an optimizer for solving this problem~(illustrated in the right of \cref{overview1}). 
In \cref{sec3:instruction_set}, we introduce how to establish a high-quality instruction set that comprises expert-level knowledge about solving optimization problems. Based on the proposed instruction set, we design a novel two-phase instruction tuning strategy to smoothly boost the performance, which is detailed in \cref{sec3:instruction_tuning}. 

\subsection{Construction of Instruction Set}\label{sec3:instruction_set}
\textbf{Task synthesis.} An optimization problem can be mathematically formulated as follows:
\begin{align}\label{eq:1}
    Minimize:& \quad f(x), \quad x=(x_1,x_2,...,x_D) \notag\\
    s.t.:&\quad g_i(x) \leq 0, \quad i=1,...,M_g& \notag\\
    &\quad h_j(x)=0,\quad j=1,...,M_h&
\end{align}
where $f(\cdot)$ is the objective function, $x$ is a $D$-dimensional vector denoting a solution, $g_i(\cdot)$ and $h_j(\cdot)$ denote $M_g$ inequality constraints and $M_h$ equality constraints respectively. Without loss of generality, we assume a minimization problem where the optimal solution $x^*$ attains the minimum objective value, adhering to all specified constraints. 

The main concern in this context is how to create an adequate amount of problem instances that possess both high quality and diversity, which is crucial for instruction tuning~\cite{formatteddata,zhou2023lima}. Since it is not feasible to gather all types of optimization problems that arise in realistic scenarios, we opt for a more feasible approach by generating synthetic problem instances. Specifically, we collect a basic function set $\digamma$ comprising many different optimization problems and a basic constraint set $\Omega$ comprising a wide variety of constraints from the well-known optimization benchmarks~\cite{boyd2004convex,cec2017cons,ma2023metabox}. Following the methodology of Mohamed et al.~(\citeyear{cec2021TR}), we synthesize a new objective function based on $K$ basic functions in $\digamma$ through two different paradigms as given by \cref{eq:2}. 1)~\textit{Composition}: this involves a linear combination of the $K$ basic functions on the complete decision space, where each $w_i$ is uniformly sampled in $[0,1]$. 2)~ \textit{Hybrid}: 
we randomly decompose $x$ into $K$ segments $s_1$ to $s_K$. The $K$ basic functions then operate on these $K$ segments, respectively, and the final objective function is the summation of these basic functions' values on the corresponding decision subspace. 
\begin{align}\label{eq:2}
    Composition: &\quad  f(x) = \sum_{i=1}^{K} w_i \cdot f_i(x) \notag \\ 
    Hybrid: & \quad f(x) = \sum_{i=1}^{K}f_i(x[s_i])
\end{align}
More concretely, we obtain a problem instance by three steps: 1) Firstly, we indicate the problem dimension $D$, the search bounds for each dimension~(e.g., $-10\leq x_i \leq10$), 
and the number of basic functions $K$; 2) Secondly, if $K=1$, we randomly select a basic function in $\digamma$ as $f(x)$, otherwise, we apply \textit{Composition}/\textit{Hybrid} paradigm to synthesize $f(x)$; and 3) Lastly, we randomly select a group of constraints $\{\{g_i\},\{h_j\}\}$ in $\Omega$. Note that step 3) is optional, as some optimization problems may not have constraints. 

In this work, we generate $3$k problem instances without constraints, denoted as $P_{\rm{nc}}$, and another $3$k problem instances with constraints, denoted as $P_{\rm{c}}$. The complete set $P$ is the union of $P_{\rm{nc}}$ and $P_{\rm{c}}$, consisting of $6$k instances. These instances showcase different characteristics of global landscapes, including unimodal or multimodal, separable or nonseparable, and symmetrical or asymmetrical. They also exhibit various local landscape properties, such as distinct properties around different local optima, continuous everywhere yet differentiable nowhere, and optima situated in flattened areas. This guarantees that the generated instances comprehensively mirror various realistic problems.   

\textbf{Knowledge gathering.} In our study, the term `knowledge' refers to expertise on how to deal with an optimization problem, which involves identifying a well-performing optimizer and configuring its hyper-parameters. After synthesizing the task set, we conduct exhaustive benchmarking to determine one effective optimizer for each instance $p \in P$. Concretely, we filter a wide range of optimizers from the published literature~\cite{optimizer-survey1,optimizer-survey2},  competitions~\cite{cec2017cons,cec2021TR,nips-challenge2020}, and benchmarks~\cite{bayesmark,ma2023metabox}. The $23$ selected optimizers form an algorithm pool, which covers various algorithm families, including Evolutionary Algorithms~(e.g., GA~\cite{ga,clune2008natural,wang2023improved}, DE~\cite{de,xu2020helper,biswas2021improving,ye2023differential}, PSO~\cite{pso,gong2015genetic,wu2022employing,lu2023double} and ES~\cite{es,ros2008simple,hansen2009benchmarking,new-es}), Bayesian Optimization~\cite{bo,wang2020learning}, Local Search strategies~\cite{kirkpatrick1983optimization,van1987simulated,xiang1997generalized,new-sa}, and Numerical Optimization methods~\cite{slsqp,conn2000trust,cobyla,new-BFGS}. 
To determine the most effective optimizer among our algorithm pool for each instance $p$, we employ a two-step process. Firstly, we perform a grid search to identify the best configuration for each optimizer on $p$ (conducted multiple times to reduce the impact of variance). Subsequently, we select the optimizer that yields the best performance among all the configured optimizers. The selected optimizer and its configuration are implemented as a piece of Python code, serving as the knowledge of the desired optimizer's implementation for instance $p$. Refer to \cref{appx:benchmark} for details of those selected optimizers~(configurations, implementations etc.) and the benchmarking process.

\begin{figure}[t]
\vskip 0.1in
\begin{center}
\centerline{\includegraphics[width=0.98\columnwidth]{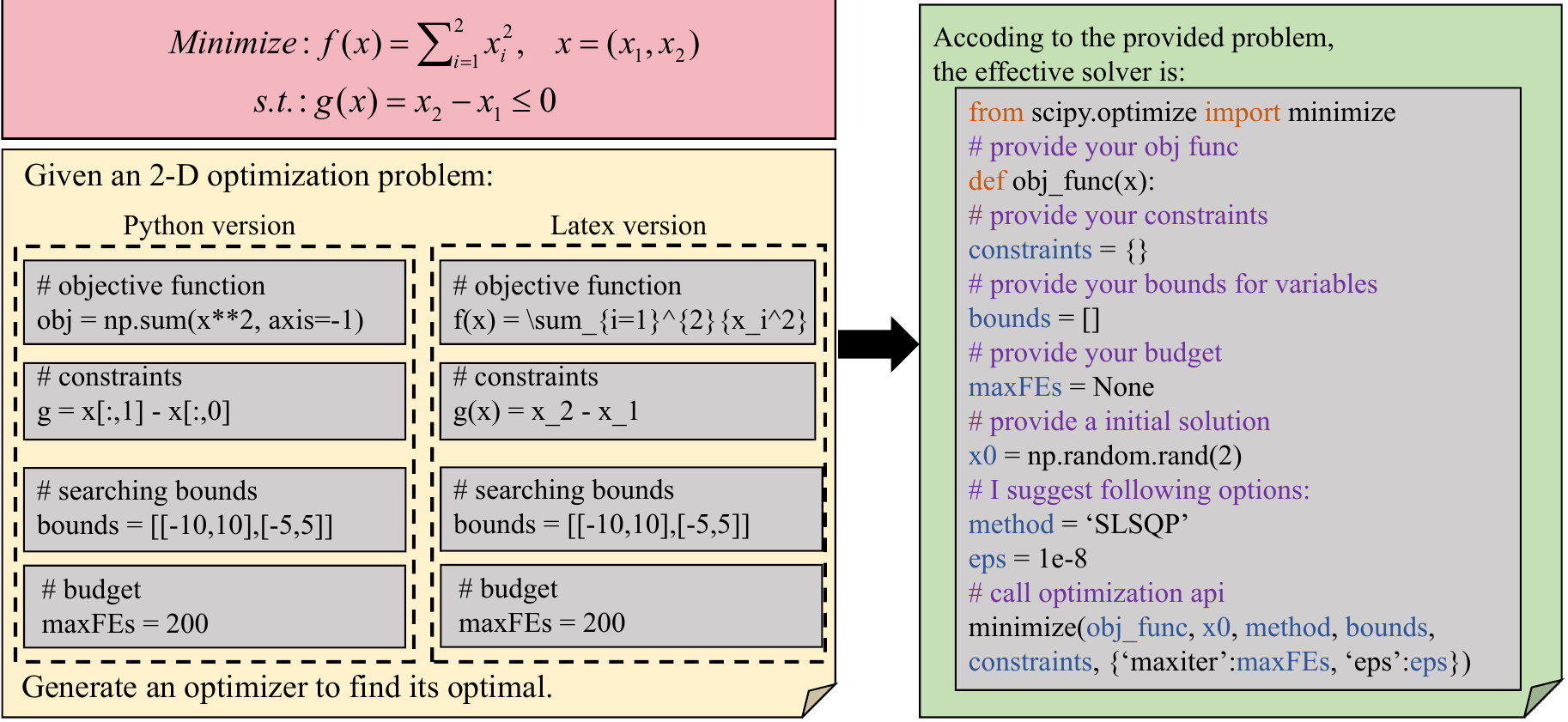}}
\vskip -0.02in
\caption{Input-output example of our instruction set. For a given problem~(red box), diverse task descriptions in Python/LaTeX formats construct the prompt~(yellow box). The code implementation of an effective optimizer is provided as the answer~(green box).}
\vspace{-5mm}
\label{fig:prompt}
\end{center}
\end{figure}

\textbf{Enhancement with diverse task descriptions.} Recent studies suggest that enhancing the diversity of the task descriptions for each task instance can lead to additional generalization gains for the instruction-tuned LLMs~\cite{formatteddata,IT2,IT-scaling}. We hence augment each problem instance in $P$ by rephrasing the writing style of its objective function and constraints. Specifically, we invited a total of $500$ university students majoring in computer science to write Python or LaTeX codes for describing a variety of optimization problems. 
Different writing patterns are observed during this process. Based on these different patterns, for each problem instance $p\in P$, we can obtain a number of rephrased versions for describing its objective function and constraints in either Python or LaTeX code. Refer to \cref{appx:rephrase} for the detailed rephrasing process and the different writing styles we have found.

After the data augmentation, we obtain the final instruction set by transforming each instance $p$, along with its rephrased versions, into a text prompt $q$~(input), and setting the source code of the selected optimizer with configurations as the answer $a$~(output). This results in an instruction set $\mathbb{I}$ comprising $32570$ pairs of input-output examples $(q,a)$, where an input-output example is illustrated in \cref{fig:prompt}.  

\subsection{Two-Phase Instruction Tuning}\label{sec3:instruction_tuning}
\textbf{Contrastive warm-up.} A key observation during the instruction set construction process in \cref{sec3:instruction_set} is that: even two prompts $q_m$ and $q_n$ are very different to each other~(e.g., they adopt different descriptions of the same problem), they can share the same desired optimizer $a$. On the contrary, for two prompts hold similar descriptions, the selected optimizers may differ. This phenomenon challenges the convergence of the models during fine-tuning. 
An appealing approach to alleviate this issue is to adopt contrastive learning to align the latent space representation for different prompts that share the same semantics. Such contrastive learning task has shown its effectiveness in several code understanding scenarios~\cite{guo2022unixcoder}. 
In LLaMoCo, we adopt constrastive learning~\cite{cl} to warm up the LLMs before instruction tuning. 

\begin{figure}[t]
\vskip 0.1in
\begin{center}
\centerline{\includegraphics[width=0.8\columnwidth]{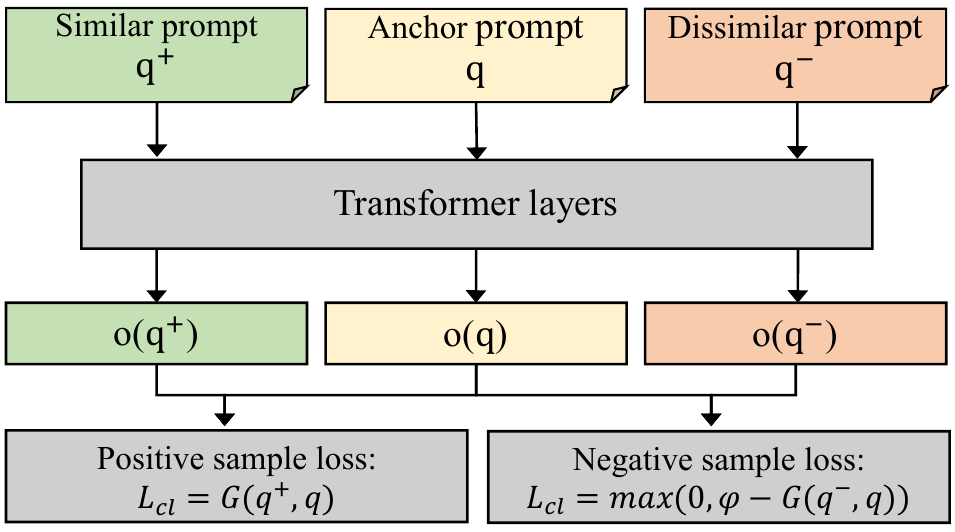}}
\caption{The workflow of contrastive warm-up strategy. Given an anchor problem prompt, we collect its similar prompts and dissimilar prompts within a mini-batch. The positive and negative sample losses are calculated on the latent space representations.}
\label{fig:contrastive_concept}
\end{center}
\vskip -0.3in
\end{figure}

The workflow of the loss calculation is illustrated in \cref{fig:contrastive_concept}. Given an anchor prompt, we collect its similar
prompts and dissimilar prompts within a mini-batch, which are then applied to calculate the positive and negative sample loss, respectively. Specifically, for the decoder-only LLMs adopted for code generation tasks in this paper, we activate the Transformer layers~\cite{vaswani2017attention} and regard the output embedding of the final self-attention block as latent space representation for the prompt $q$. 

In LLaMoCo, we measure the distance between two prompts $q_m$ and $q_n$, denoted as $G(q_m, q_n)$, by considering the cosine similarity between their latent space representations $\overrightarrow{o}(q_m)$ and $\overrightarrow{o}(q_n)$:
\begin{equation}
    G(q_m, q_n) = \frac{1}{2} \left(1-\frac{\vec{o}\left(q_m\right) \cdot \vec{o}\left(q_n\right)}{\left\|\vec{o}\left(q_m\right)\right\|\left\|\vec{o}\left(q_n\right)\right\|}\right)
\end{equation}
The above distance $G(q_m, q_n) \in [0,1]$. Then, the contrastive loss of $q_m$ and $q_n$, denoted as $L_{\rm{cl}}(q_m,q_n)$, is as
\begin{equation}
L_{\rm{cl}}= 
\begin{cases}
G(q_m,q_n)& a_m = a_n\\
\max (0, \varphi-G(q_m,q_n))& a_m \neq a_n
\end{cases}
\end{equation}
where $a_m$ and $a_n$ are the corresponding designated optimizer of $q_m$ and $q_n$, respectively, $\varphi$ is a margin parameter. By minimizing $L_{\rm{cl}}$, we could efficiently pull together the representations of two prompts which share the same desired optimizer yet have different forms, and vice versa. In LLaMoCo, we consume a small number of epochs to warm up the fine-tuning of LLMs by $L_{\rm{cl}}$ and then instruction-tune the LLMs with the normal language modelling loss for next-token prediction~\cite{wolf2019huggingface}. We note that the contrastive warm-up phase does not require context generation, hence the time cost is relatively smaller compared with the subsequent instruction tuning phase. We validate the effectiveness of this contrastive learning phase in \cref{sec:ablation}.

\textbf{Balanced data sampling.} The instruction set $\mathbb{I}$ exhibits certain imbalance in the distribution of data. Notably, we observe that several optimizers dominate on thousands of problem instances, while the others only outperform on a few problem instances. Dealing with imbalanced data poses a challenge during the process of fine-tuning models~\cite{balance-kdd,llm-survey}. 
To address the issue, we follow the example-proportional mixing strategy~\cite{distribution1} to re-balance the data distribution in $\mathbb{I}$. Each data pair $(q,a)$ is sampled with a probability $\rho$ as:
\begin{equation}
    \rho(q,a) = \frac{1}{N_a\times N_{q,a}}
\end{equation}
where $N_a$ denotes the number of optimizers in the gathered algorithm pool, $N_{q,a}$ denotes the number of instances whose desired optimizer is $a$. In this way, the number of sampled pairs dominated by each optimizer is approximately equal in each training epoch. Note that we apply this strategy in both the contrastive warm-up phase and the instruction tuning phase. The approach aids in avoiding biased training of the LLMs and enables them to effectively learn the knowledge from minority instances. 
In addition, a homogeneous mini-batch sampling strategy is applied, due to the space limitation, it is presented in \cref{appx:training_detail}.


\section{Results and Discussions}\label{sec4:exp}
\subsection{Experimental Setup}
\textbf{Fundamental models.} We adopt CodeGen-Mono~($350$M), Phi-$2$~($2.7$B) and Code Llama~($7$B) as fundamental models and fine-tune them on our instruction set. The reasons are two-fold: 1) these models show robust programming language reasoning and code generation ability, serving as a good start point for the code-to-code scenario in our work; 2) the relatively small model size helps to reduce computational resources required for training and deploying.

\textbf{Training settings.} For generating the task set $P$, the problem dimension $D$ for each $p_i$ is randomly chosen from $[2,50]$, and the number of components $K$ is randomly chosen from $[1,5]$. We randomly split the instruction set $\mathbb{I}$ into a training set $\mathbb{I}_{\rm{train}}$ with $30$k input-output pairs and a test set $\mathbb{I}_{\rm{eval}}$ with the rest examples. For our two-phase instruction tuning, we deploy $5$ epochs of contrastive warm-up and $20$ epochs of instruction tuning for all fundamental models. Specifically,
we first apply \textit{SGD}~\cite{sgd} with a fixed learning rate $5\times 10^{-4}$ in the contrastive warm-up phase, alongside $\varphi = 0.3$. Then, we apply \textit{AdamW}~\cite{adamw} to optimize the LLMs in the instruction tuning phase. During the initial $1$k iterations, the learning rate gradually increases from $0$ to $5\times 10^{-4}$ in a linear manner. Subsequently, it decreases to $0$ according to a cosine schedule. The batch size in both phases is set to $4$. Note that we fine-tune the CodeGen-Mono~($350$M) with full parameters, but apply LoRA~\cite{hu2022lora} to fine-tune the larger Phi-$2$~($2.7$B) and Code Llama~($7$B) models, with the rank $r=8$, scaling factor $\alpha = 32$, and a dropout rate of $0.05$. All experiments are performed on a platform with an Intel(R) Xeon(R) Gold $6348$ CPU, $504$GB RAM and a Nvidia A$800$~($80$GB) GPU. Upon the settings, the training duration for CodeGen is one day, whereas Phi-$2$ and Code Llama require $2.5$ days and $4$ days of training, respectively.

\textbf{Competitors.} We include two solution-to-solution approaches, OPRO~\cite{opro} and LMEA~\cite{opro-ntu}, which prompt pre-trained LLMs~(e.g., GPT-$4$ Turbo) repeatedly to generate and improve solutions for the given problems. 
Compared to OPRO, LMEA additionally engineered its prompt with an explicit indication of using some evolutionary operators to let LLMs act as an evolutionary optimizer for performance boost. We also include three general LLMs for code generation, namely Code Llama-$7$B~\cite{codellmma}, Llama $2$-$70$B~\cite{llama2}, and GPT-$4$ Turbo~\cite{gpt-4t}. We prompt these three general LLMs with the same format as in our instruction set $\mathbb{I}$ to generate an optimizer for each problem instance. The configurations of the competitors are set by default according to the corresponding references. 

\begin{table*}[t]
\centering
\caption{Results of different approaches in terms of \textbf{Code Error Rate~(Err.)}, \textbf{Code Recovery Cost~(Rec.)}, \textbf{Optimization Performance~(Perf.)}, and \textbf{Computational Overhead~(Comp.)} on the unconstrained problems ($\mathbb{I}_{\rm{eval}}/P_{\rm{c}}$), constrained problems ($\mathbb{I}_{\rm{eval}}/P_{\rm{nc}}$), and all test problems ($\mathbb{I}_{\rm{eval}}$), where ``-'' denotes that the approach does not generate code (it follows a solution-to-solution paradigm). 
}\vskip 0.2in
\label{tab:overall}
\resizebox{0.9\textwidth}{!}{%
\begin{tabular}{c|c|cc|ccc|ccc}
 \hline
  \multirow{3}{*}{Testset} &
  \multirow{3}{*}{Metrics} &
  \multicolumn{2}{c|}{\multirow{2}{*}{Prompt for Solution}} &
  \multicolumn{3}{c|}{\multirow{2}{*}{Prompt for Optimizer}} &
  \multicolumn{3}{c}{\multirow{2}{*}{Our LLaMoCo}} \\
 &
   &
  \multicolumn{2}{c|}{} &
  \multicolumn{3}{c|}{} &
  \multicolumn{3}{c}{} \\ \cline{3-10}
 &
   &
OPRO &
LMEA &
GPT-$4$ Turbo &
Code Llama-$7$B &
Llama2-$70$B &
\begin{tabular}[c]{@{}c@{}}LLaMoCo-S\end{tabular} &
\begin{tabular}[c]{@{}c@{}}LLaMoCo-M\end{tabular} &
\begin{tabular}[c]{@{}c@{}}LLaMoCo-L\end{tabular} 
\\ \hline
\multirow{4}{*}{\begin{tabular}[c]{@{}c@{}}$\mathbb{I}_{\rm{eval}}/P_{\rm{c}}$\end{tabular}}
 &
  \begin{tabular}[c]{@{}c@{}}Err. $\downarrow$\end{tabular}
   & -
   & -
   & 43.333\%
   & 98.184\%
   & 99.673\%
   & 5.437\%
   & \textbf{4.414\%}
   & 4.697\%
   \\ \cline{2-2}
 &
  \begin{tabular}[c]{@{}c@{}}Rec. $\downarrow$\end{tabular}
   & -
   & -
   & 9.942\%
   & 67.857\%
   & 62.232\%
   & \textbf{9.684\%}
   & 10.101\%
   & 9.947\%
   \\ \cline{2-2}
   &
  \begin{tabular}[c]{@{}c@{}}Perf. $\uparrow$\end{tabular}
   & 29.499\% 
   & 20.350\% 
   & 71.783\% 
   & 14.089\% 
   & 18.922\% 
   & 85.360\% 
   & \textbf{86.412\%} 
   & 85.810\% 
   \\ \cline{2-2}
   &
  \begin{tabular}[c]{@{}c@{}}Comp. $\downarrow$\end{tabular}
   & 115k
   & 249k
   & 3.4k
   & 1.7k
   & \textbf{1.5k}
   & 2.3k
   & 2.3k
   & 2.3k
   \\ \hline
\multirow{4}{*}{\begin{tabular}[c]{@{}c@{}}$\mathbb{I}_{\rm{eval}}/P_{\rm{nc}}$\end{tabular}}
 &
  \begin{tabular}[c]{@{}c@{}}Err. $\downarrow$\end{tabular}
   & -
   & -
   & 39.944\%
   & 90.474\%
   & 99.521\%
   & \textbf{5.697\%}
   & 6.130\%
   & 5.977\%
   \\ \cline{2-2}
 &
  \begin{tabular}[c]{@{}c@{}}Rec. $\downarrow$\end{tabular}
   & -
   & -
   & 16.463\%
   & 44.938\%
   & 49.202\%
   & 11.861\%
   & \textbf{10.443\%}
   & 10.584\%
   \\ \cline{2-2}
   &
  \begin{tabular}[c]{@{}c@{}}Perf. $\uparrow$\end{tabular}
   & 4.514\% 
   & 7.541\% 
   & 75.678\% 
   & 46.968\% 
   & 22.460\% 
   & 77.576\% 
   & 79.718\% 
   & \textbf{83.404\%} 
   \\ \cline{2-2}
   &
  \begin{tabular}[c]{@{}c@{}}Comp. $\downarrow$\end{tabular}
   & 115k
   & 249k
   & 3.5k
   & \textbf{2.0k}
   & \textbf{2.0k}
   & 2.5k
   & 2.5k
   & 2.5k
 \\ \hline
\multirow{4}{*}{\begin{tabular}[c]{@{}c@{}}$\mathbb{I}_{\rm{eval}}$\end{tabular}}
 &
  \begin{tabular}[c]{@{}c@{}}Err. $\downarrow$\end{tabular} 
   & -
   & -
   & 41.667\%
   & 95.156\%
   & 99.617\%
   & 5.580\%
   & \textbf{5.434\%}
   & 5.509\%
   \\ \cline{2-2}
 &
  \begin{tabular}[c]{@{}c@{}}Rec. $\downarrow$\end{tabular} 
   & -
   & -
   & 13.072\%
   & 57.001\%
   & 55.717\%
   & 10.826\%
   & \textbf{10.349\%}
   & 10.461\%
   \\ \cline{2-2}
   &
  \begin{tabular}[c]{@{}c@{}}Perf. $\uparrow$\end{tabular}
   & 17.821\% 
   & 14.762\% 
   & 74.248\% 
   & 29.717\% 
   & 20.579\% 
   & 81.843\% 
   & 83.369\% 
   & \textbf{83.451\%} 
   \\ \cline{2-2}
   &
  \begin{tabular}[c]{@{}c@{}}Comp. $\downarrow$\end{tabular}
   & 115k
   & 249k
   & 3.5k
   & 1.9k
   & \textbf{1.7k}
   & 2.4k
   & 2.4k
   & 2.4k
\\ \hline
\end{tabular}%
}
\vspace{-3mm}
\end{table*}

\textbf{Performance metrics.} When evaluating the performance of LLMs for optimization, we consider four metrics: 1) the \textit{code error rate}, which indicates the proportion of problems for which the LLMs generate optimization codes with bugs (lower values are preferable); 2) the \textit{code recovery cost}, which measures the proportion of lines of code that need to be corrected in order to fix the bugs in the erroneous codes (lower values are preferable); 3) the average \textit{optimization performance} on the test problems (higher values are preferable); and 4) the average \textit{computational overhead} for solving a problem, which is determined by the number of tokens used for both the input and output of LLMs (lower values are preferable). These four metrics could provide a comprehensive evaluation on existing baselines and our LLaMoCo in aspects of code generation robustness, optimization performance and runtime complexity. The detailed calculations for these metrics can be found in Appendix \ref{appx:calculation}.


\begin{table}[t]

\centering
\caption{Performance comparison on realistic problems.}
\vskip 0.15in
\label{tab:realistic}
\resizebox{0.8\columnwidth}{!}{%
\begin{tabular}{c|c|c|c}
 \hline
Metrics &
OPRO &
GPT-$4$ Turbo &
\begin{tabular}[c]{@{}c@{}}LLaMoCo-S\end{tabular}
\\ \hline
  \begin{tabular}[c]{@{}c@{}}Err. $\downarrow$\end{tabular}
   & -
   & 79.483\%
   & \textbf{4.168\%}
   \\ \cline{1-1}
  \begin{tabular}[c]{@{}c@{}}Rec. $\downarrow$\end{tabular}
   & -
   & 22.985\%
   & \textbf{7.426\%}
   \\ \cline{1-1}
  \begin{tabular}[c]{@{}c@{}}Perf. $\uparrow$\end{tabular}
   & 73.995\% 
   & 59.174\% 
   & \textbf{87.227\%} 
   \\ \cline{1-1}
  \begin{tabular}[c]{@{}c@{}}Comp. $\downarrow$\end{tabular}
   & 241k
   & 3.6k
   & \textbf{2.5k}
   \\ \hline
\end{tabular}%
}
\end{table}
\subsection{Performance Analysis}
We use LLaMoCo-S(mall), -M(edium) and -L(arge) to denote the fine-tuned CodeGen-Mono~($350$M), Phi-$2$~($2.7$B) and Code Llama~($7$B) models on $\mathbb{I}_{\rm{train}}$, respectively. 

\textbf{Performance on test sets.} First, we evaluate the performance of our fine-tuned LLMs and the competitors on three test sets, $\mathbb{I}_{\rm{eval}}/P_{\rm{c}}$, $\mathbb{I}_{\rm{eval}}/P_{\rm{nc}}$, and $\mathbb{I}_{\rm{eval}}$ that represent the unconstrained task set, constrained task set, and the complete set mixing unconstrained and constrained tasks, respectively, each with $5$ independent runs. The results in terms of the four metrics are reported in \cref{tab:overall}, which show that: 

1) The LLMs fine-tuned by our LLaMoCo framework consistently achieve superior performance, which validates that instruction tuning the general LLMs with moderate expert-level knowledge would gain substantial performance reinforcement in optimization. 
For example, LLaMoCo-L fine-tuned on the Code Llama~($7$B) demonstrate an optimization performance boost from $29.717\%$ to $81.843\%$ on~$\mathbb{I}_{\rm{eval}}$.

2) Although LLaMoCo-S is fine-tuned from a relatively small fundamental model, it achieves competitive performance to those of LLaMoCo-M and LLaMoCo-L. This may reveal a potential marginal effect in instruction tuning, since the data scale should match the capacity of the model. 

3) The solution-to-solution approaches OPRO and LMEA achieve unsatisfactory performance on our complex optimization task sets. 
Considering the tremendous tokens these approaches consume to solve one optimization problem through iteratively prompting solutions, both the efficacy and efficiency~(as shown in the `Perf.' and `Comp.' rows of \cref{tab:overall}) of them require further improvement.

4) Among the three `prompt for optimizer' models we compared, the GPT-$4$ Turbo dominates the other two, which shows the powerfulness of a general-purpose LLM with high capacity. Nevertheless, it still underperforms our domain-specific LLaMoCo. Our models effectively reduce the error rates and the required recovery efforts for generating the codes of an optimizer through the instruction tuning. 
Meanwhile, note that the Code Llama~($7$B) model achieves better overall performance than the Llama $2$~($70$B) model in our experiments. The above observations validate that, although LLMs with larger capacity may show strong performance for solving general tasks, a smaller model could be sufficient to be fine-tuned as a domain-specific task solver. 

\textbf{Zero-shot performance on realistic problems.} We introduce a realistic optimization problem set collected by Kumar et al.~(\citeyear{kumar2020test}) to further evaluate the zero-shot generalization performance of the LLMs fine-tuned by our LLaMoCo. This problem set serves as an ideal testbed for our framework for two reasons: 1) an optimizer that performs very well on synthetic benchmark suites may not provide robust performance on real-world problems, and 2) this set of problems shows very different structures compared to our synthetic problem set $P$. These problems come from various real-world scenarios including the industrial chemical process, mechanical engineering, and the power system, many of them feature high-dimensional problem spaces and complicated constraints. As an illustration, we test OPRO, GPT-$4$ Turbo and our LLaMoCo-S on these realistic problems~(integrate their problem definitions into our formatted prompts) for $5$ independent runs. The results in \cref{tab:realistic} demonstrate the best generalization performance and hence the practical availability of our LLaMoCo.

\subsection{Ablation study}\label{sec:ablation}
\begin{figure}[t]
\centering
\vskip -0.08in
\subfigure{
\label{fig:rephrase-lang}
\includegraphics[width=0.46\columnwidth]{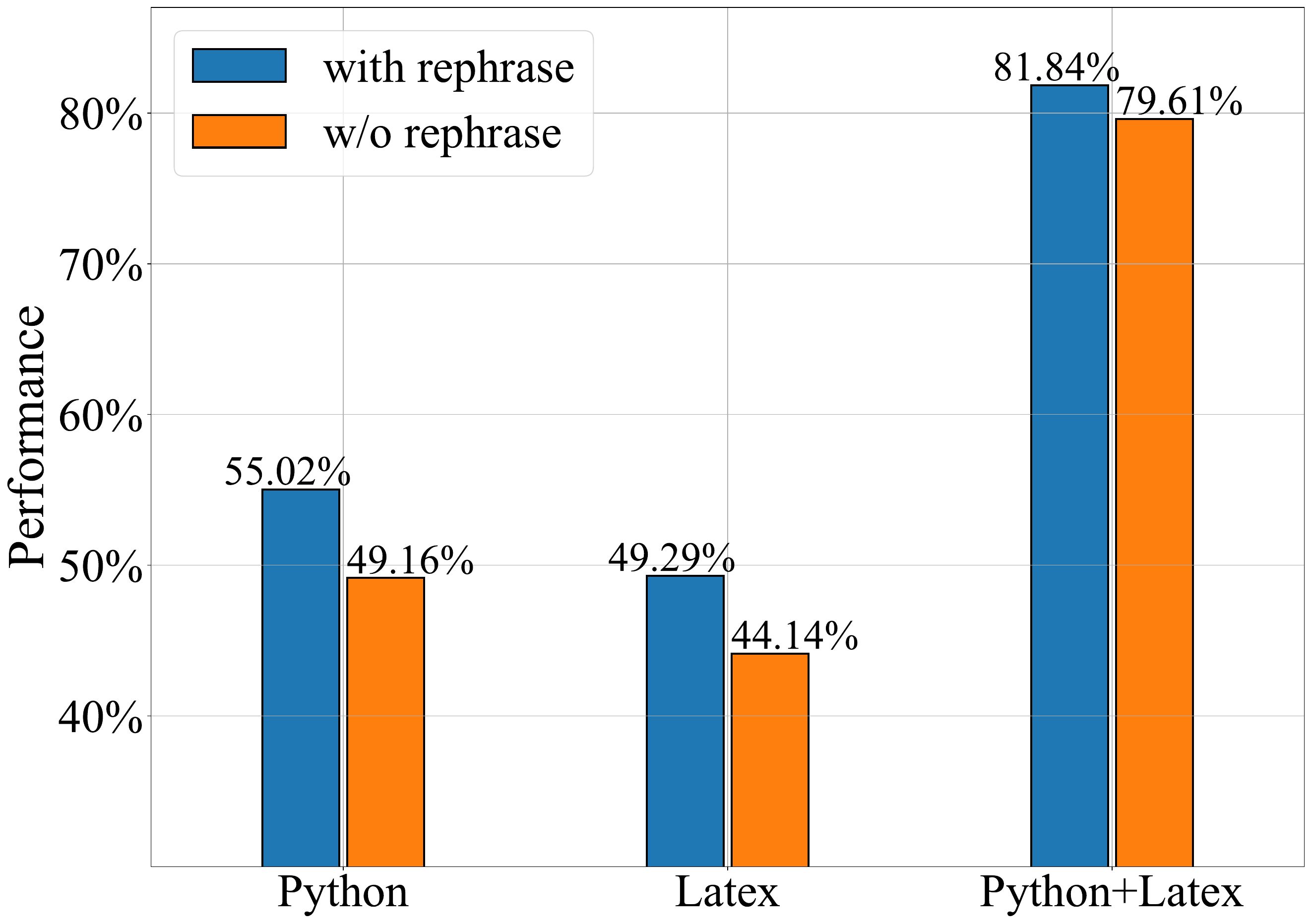}
}
\hfill
\subfigure{
\label{fig:quality}
\includegraphics[width=0.46\columnwidth]{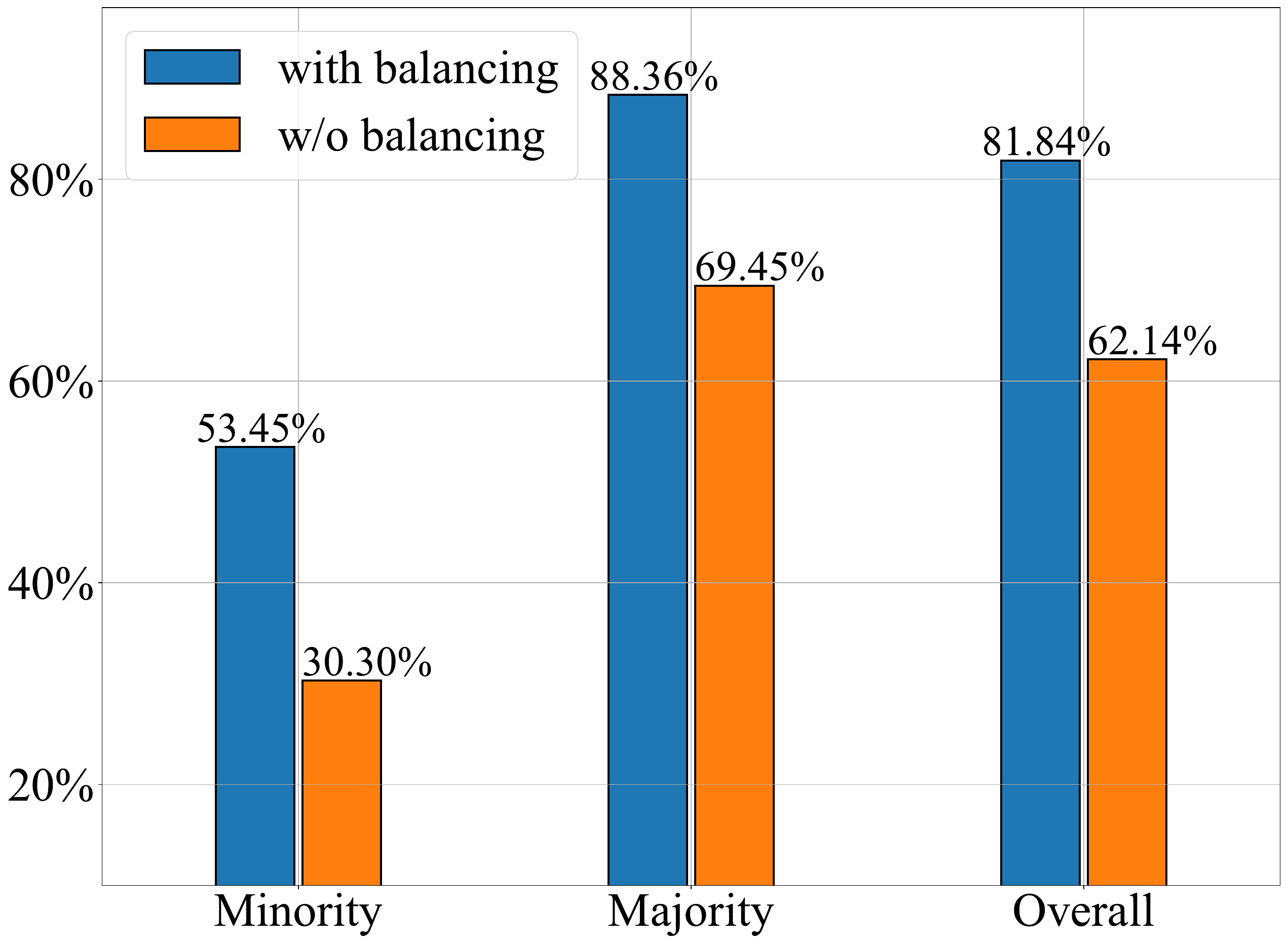}
}
\vskip -0.08in
\caption{Effectiveness of the diversity enhancement strategy (\textbf{Left}) and the data distribution balancing strategy (\textbf{Right}).
}
\label{fig:diver-quali}
\end{figure}

\textbf{Diversity enhancement.} To improve the generalization of the fine-tuned LLMs in LLaMoCo, we enrich the task descriptions for each problem instance by augmenting the description of its objective function and constraints with Python or LaTeX codes of different writing styles. We illustrate the effect of this procedure in the left of \cref{fig:diver-quali} by showing the optimization performance of six LLaMoCo-S models trained on pure Python, pure LaTeX and Python+LaTeX data, with or without the diversity enhancement by rephrasing. The results show that providing multi-lingual descriptions of optimization problems significantly boosts the generalization performance, while rephrasing each description with multiple writing styles further enhances the final training results.

\begin{figure}[t]
\centering
\subfigure{
\label{fig:contrastive-loss}
\includegraphics[width=0.43\columnwidth]{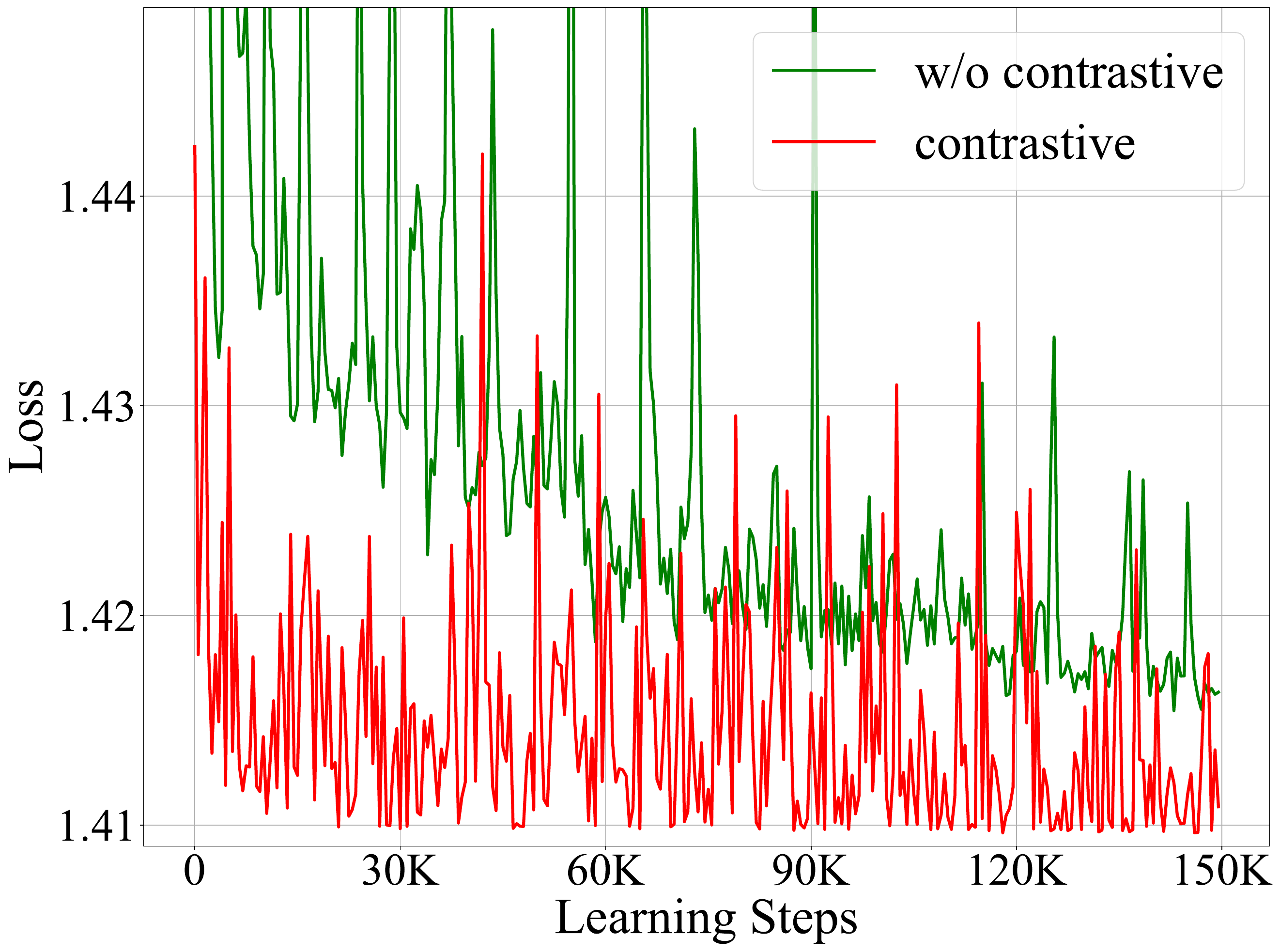}
}
\subfigure{
\label{fig:contrastive-perf}
\includegraphics[width=0.43\columnwidth]{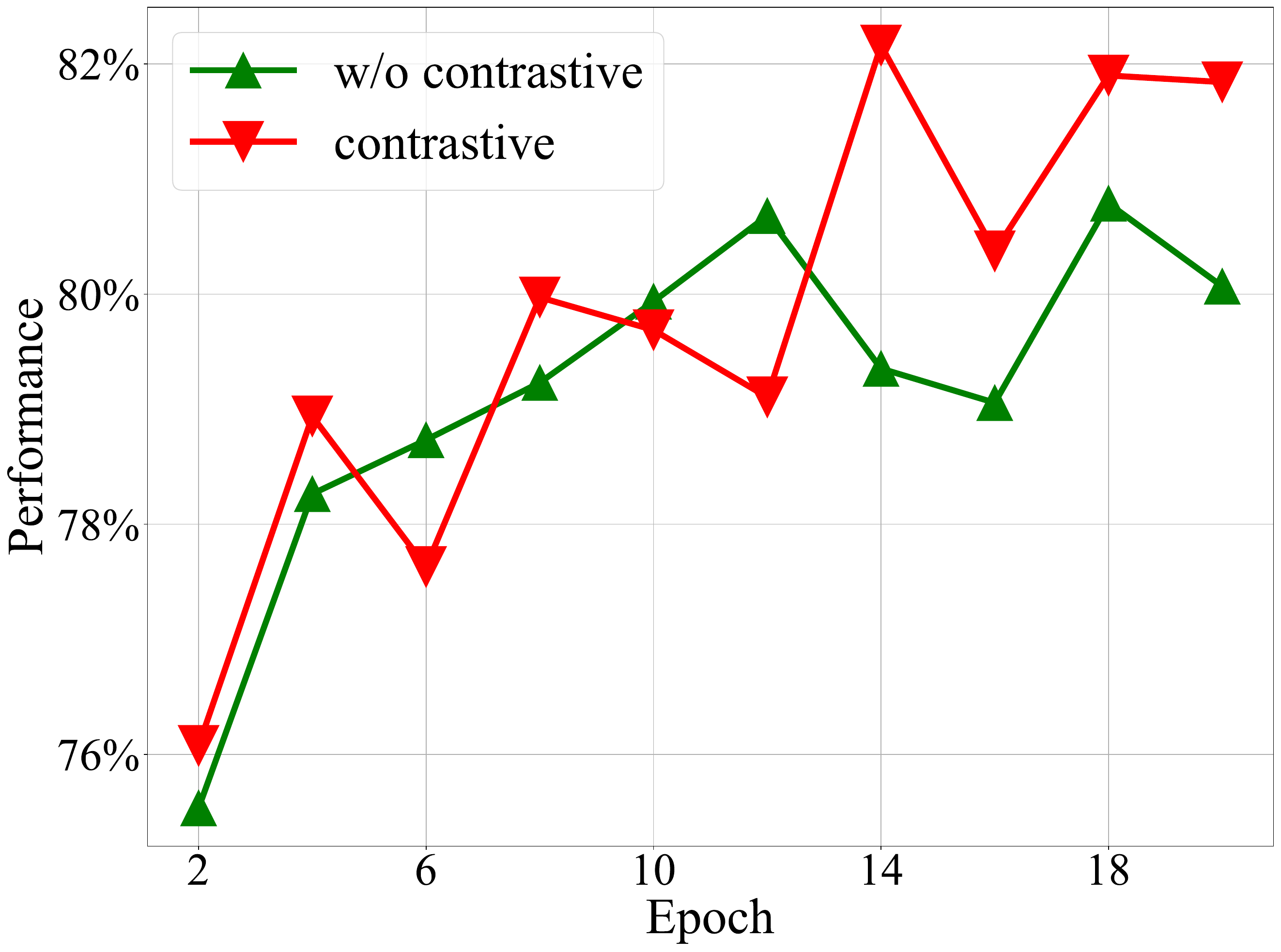}
}
\vskip -0.08in
\caption{Effectiveness of the contrastive warm-up strategy on training curves (\textbf{Left}) and performance gains (\textbf{Right}).}
\label{fig:contrastive}
\end{figure}

\textbf{Contrastive warm-up.}
The contrastive warm-up phase in our proposed two-phase instruction tuning strategy~(see \cref{sec3:instruction_tuning}) aims to reduce the cross-modal ambiguity by aligning the latent space representations of different prompts that share the same desired optimizer. We illustrate the training curves and performance gain curves on $\mathbb{I}_{\rm{eval}}$ with or without the contrastive warm-up during the instruction tuning phase in \cref{fig:contrastive}, where LLaMoCo-S is applied as a showcase. The results show that incorporating such a contrastive warm-up strategy aids in accelerating the convergence of the subsequent instruction tuning. 
Furthermore, it is advantageous for the LLMs to generate accurate codes and enhance the overall optimization performance.

\textbf{Balanced data sampling.} In LLaMoCo, we address the imbalanced data distribution (caused by dominate optimizers) through performing example-proportional sampling on $\mathbb{I}_{\rm{train}}$. To investigate its effectiveness, we train two LLaMoCo-S models on $\mathbb{I}_{\rm{train}}$, with or without the data balancing strategy, respectively. The optimization performance of the two models is presented in the right of \cref{fig:diver-quali}, by separately considering the majority instances (which request the dominating optimizers), the minority instances (which request the others), and the overall instances of $\mathbb{I}_{\rm{eval}}$. The results consistently show that keeping a balanced training data distribution significantly boosts performance.


\subsection{Open-Ended Discussion: Is GPT-4 a True Optimization Expert?} 
Considering the competitive performance of GPT-$4$ on optimization tasks, as shown in \cref{tab:overall}, we delve into whether GPT-$4$ can be deemed as a genuine optimization expert. Upon viewing the optimization codes generated by GPT-$4$ for both the test and the realistic problem set, a noteworthy pattern emerges. GPT-$4$ consistently leans towards generating a specific numerical optimizer, SLSQP~\cite{slsqp}, for almost all tested problems. While SLSQP is a classical solver for convex quadratic programming and is included in our chosen advanced optimizers, our benchmarking results identify that on a proportion of tested problems, it underperforms the others such as the  Vanilla DE~\cite{de}. To investigate further, we experiment by providing GPT-$4$ with a hint to use Vanilla DE to solve these specific problems. Surprisingly, GPT-$4$ successfully outputs a code implementation of DE and achieves competitive results. This observation suggests that while GPT-$4$ may have included sufficient domain knowledge on how to solve optimization problems, it still exhibits an underfitting issue concerning how to solve a `particular' problem. This underscores the importance of the LLaMoCo framework for fine-tuning general LLMs to fit the task of generating an appropriate optimizer tailored for specific problem instances.

\section{Conclusion}\label{sec5:conclusion}
We introduce LLaMoCo, the first instruction-tuning framework to adapt general LLMs to function as expert-level systems to solve optimization problems. To achieve this, we meticulously construct an instruction set with more than $30$k demonstration examples and then employ a novel two-phase instruction tuning strategy to fine-tune a series of LLMs. The results show that our models consistently outperform existing approaches. Notably, we observe that a relatively small LLM is sufficient to be tuned as an expert-level optimization code generator superior to GPT-$4$. As a preliminary exploratory research endeavour, LLaMoCo certainly has limitations, such as the need to augment the instruction set with more instances to enhance generalization performance. Additionally, we consider enhancing the LLMs fine-tuned by LLaMoCo through further alignment tuning as a promising future direction.
\section*{Impact Statements}
This paper presents work whose goal is to advance the field of Machine Learning. There are many potential societal consequences of our work, none of which we feel must be specifically highlighted here.
\bibliography{example_paper}
\bibliographystyle{icml2024}

\newpage
\appendix
\onecolumn
\section{Benchmarking for Knowledge Gathering}\label{appx:benchmark}
\subsection{Optimizer Pool and The Used Assets}
To match each problem instance in the generated problem set $P$ with an appropriate optimizer with corresponding code implementation, we construct an optimizer pool $\Lambda$ which integrates 23 well-performing optimizers from various algorithm families. These selected optimizers can be divided into two groups considering their compatibility for constraint handling. We briefly list the two groups as below:

\textbf{Unconstrained group $\Lambda_{uc}$:} Simulated Annealing~\cite{kirkpatrick1983optimization}, Vanilla PSO~\cite{pso}, Vanilla DE~\cite{de}, Dual Annealing~\cite{xiang1997generalized}, SAMR-GA~\cite{clune2008natural}, SEP-CMA-ES~\cite{ros2008simple}, BIPOP-CMA-ES~\cite{hansen2009benchmarking}, DEAP-DE~\cite{fortin2012deap}, Vanilla BO~\cite{bo}, GLPSO~\cite{gong2015genetic}, MMES~\cite{new-es}, LA-MCTS~\cite{wang2020learning}, MadDE~\cite{biswas2021improving}, sDMS-PSO~\cite{wu2022employing}, AMCDE~\cite{ye2023differential}, NSA~\cite{new-sa}.

\textbf{Constrained group $\Lambda_{c}$:} SLSQP~\cite{slsqp}, Trust-Constr~\cite{conn2000trust},
COBYLA~\cite{cobyla}, L-BFGS-B~\cite{morales2011remark}, HECO-DE~\cite{xu2020helper}, DTPSO~\cite{lu2023double}, GA-TDX~\cite{wang2023improved}.

We benefit from open-source libraries, including DEAP~\cite{fortin2012deap}, PyPop7~\cite{duan2022pypop7}, evosax~\cite{lange2023evosax}, SciPy~\cite{virtanen2020scipy} and Scikit-Optimizer~\cite{bo2014bo} etc., for the easy implementation of the selected optimizers. We list the codebases we adopt for the implementation of these optimizers and their licenses in \cref{tab:licenses}. We note that the development and deployment of our framework strictly follow those licenses.
\begin{table}[h]
    \caption{Used assets and their licenses}
    \label{tab:licenses}
    \centering
    \resizebox{0.9\textwidth}{!}{
    \begin{tabular}{ccc}
    \\
    \toprule
         Asset & Codebase & License \\
         \midrule
         DEAP-DE~\cite{fortin2012deap} & \multirow{2}{*}{DEAP~\cite{fortin2012deap}} & \multirow{2}{*}{LGPL-3.0 License}\\
         Vanilla PSO~\cite{pso} &  & \\
         \midrule
         SAMR-GA~\cite{clune2008natural} & \multirow{3}{*}{evosax~\cite{lange2023evosax}} & \multirow{3}{*}{Apache-2.0 license}\\
         BIPOP-CMA-ES~\cite{hansen2009benchmarking} &  & \\
         Simulated Annealing~\cite{kirkpatrick1983optimization} &  & \\
         \midrule
         SEP-CMA-ES~\cite{ros2008simple} & \multirow{4}{*}{PyPop7~\cite{duan2022pypop7}} & \multirow{2}{*}{GPL-3.0 license}\\
         MMES~\cite{new-es} &  & \\
         LA-MCTS~\cite{wang2020learning} &  & \\
         NSA~\cite{new-sa} &  & \\
         \midrule
         Dual Annealing~\cite{xiang1997generalized} & \multirow{3}{*}{SciPy~\cite{virtanen2020scipy}} & \multirow{3}{*}{BSD-3-Clause license}\\
         SLSQP~\cite{slsqp} &  & \\
         COBYLA~\cite{cobyla} &  & \\
         \midrule
         Vanilla BO~\cite{bo} & Scikit-Optimizer~\cite{bo2014bo} & BSD-3-Clause License\\
         \midrule
         GA-TDX~\cite{wang2023improved} & \multirow{8}{*}{\begin{tabular}[c]{@{}c@{}}advanced-global-optimizers\\\url{https://pypi.org/project/advanced-global-optimizers/}\end{tabular}} & \multirow{8}{*}{MIT License}\\
         Vanilla DE~\cite{de} &  & \\
         MadDE~\cite{biswas2021improving} &  & \\
         AMCDE~\cite{ye2023differential} &  & \\
         HECO-DE~\cite{xu2020helper} &  & \\
         GLPSO~\cite{gong2015genetic} &  & \\
         sDMS-PSO~\cite{wu2022employing} &  & \\
         DTPSO~\cite{lu2023double} &  & \\
         \bottomrule
    \end{tabular}}
\end{table}

\subsection{Benchmarking Process}
The benchmarking process aims to find an appropriate configured optimizer for each problem instance $p \in P$. To this end, for each optimizer $\Lambda_k \in \Lambda$, we span a grid-search configuration space $C_k$ based on its tunable hyper-parameters, which is listed in \cref{tab:optimizer}. Take DEAP-DE~\cite{fortin2012deap} as an example, it has three hyper-parameters and each of them has 5 optional values~(pre-defined by us). We hence span the $C_k$ of DEAP-DE as a configuration space comprising $5 \times 5 \times 5 = 125$ configurations, each denoted as $C^j_k$. We now establish the target of our benchmarking process: 
\begin{equation}
    \mathop{\arg\max}\limits_{\Lambda_k \in \Lambda,C^j_k \in C_k} \mathbb{E}\left[eval(p,\Lambda_k,C^j_k)\right]\nonumber
\end{equation}
where $p$ denotes the tested problem instance, $eval$ denotes the final optimization performance by calling $\Lambda_k$ with configuration $C^j_k$ to solve $p$, and $\mathbb{E}$ denotes the expectation of the optimization performance, which is unbiased-estimated by $5$ independent runs in this work. For constrained problems, we benchmark $\Lambda_{c}$, while for unconstrained problems we benchmark $\Lambda_{nc}$. We note that the benchmarking process for each problem instance may encounter execution failures, e.g., some optimizers in $\Lambda_{c}$ can not handle equality constraints, some optimizers in $\Lambda_{nc}$ are incompatible with non-convex problems, BO optimizers are extremely time-consuming on high-dimensional problems. When failures occur, the corresponding $eval(p,\Lambda_k,C^j_k)$ is set to $0$. After benchmarking $\Lambda$ on $P$, we provide a configured optimizer $a(\Lambda_k,C^j_k)$, and the corresponding code implementation as the desired optimizer for each $p$.
\scriptsize
\begin{longtable}[h]{|c|c|c|c|}
\caption{Configurations and the hyperparameter tuning settings of the optimizers.}
\label{tab:optimizer}\\
\hline
Type 
& Algorithm
& Parameters
& Search range
\\ \hline
\endfirsthead
\multicolumn{4}{c}%
{{\bfseries Table \thetable\ continued from previous page}} \\ \hline
Type & Algorithm & Parameters & Search range\\ \hline
\endhead
\multirow{9}{*}{GA}
& \multirow{5}{*}{\begin{tabular}[c]{@{}c@{}}SAMR-GA\\\cite{clune2008natural}\end{tabular}}
& NP & [10, 20, 50, 100, 200]
\\
&
& elite\_ratio & 0.0
\\
&
& sigma\_init & [0, 0.5, 1]
\\
&
& sigma\_meta & [1, 2, 3, 4, 5]
\\
&
& sigma\_best\_limit & [0.0001, 0.001, 0.1]
\\ \cline{2-4}
& \multirow{4}{*}{\begin{tabular}[c]{@{}c@{}}GA-TDX\\\cite{wang2023improved}\end{tabular}}
& beta & [0.1, 0.2, 0.3, 0.4, 0.5]
\\
&
& gamma & [1, 3, 5, 7, 9]
\\
&
& m & 1e10
\\
&
& NP & [10, 20, 50, 100, 200]
\\ \hline
\multirow{34}{*}{DE}
& \multirow{6}{*}{\begin{tabular}[c]{@{}c@{}}Vanilla DE\\\cite{de}\end{tabular}}
& NP & [10, 20, 50, 100, 200]
\\
&
& F & [0, 0.5, 0.9]
\\
&
& Cr & [0, 0.5, 0.9]
\\
&
& mutation & \begin{tabular}[c]{@{}c@{}}\{best1, best2, rand2, current2rand,\\ current2best, rand2best2\}\end{tabular}
\\
&
& bound & \{clip, periodic, reflect, rand\}
\\ \cline{2-4}
& \multirow{3}{*}{\begin{tabular}[c]{@{}c@{}}DEAP-DE\\\cite{fortin2012deap}\end{tabular}}
& NP & [10, 20, 50, 100, 200]
\\
&
& F & [0.1, 0.3, 0.5, 0.7, 0.9]
\\
&
& Cr & [0.1, 0.3, 0.5, 0.7, 0.9]
\\ \cline{2-4}
& \multirow{8}{*}{\begin{tabular}[c]{@{}c@{}}HECO-DE\\\cite{xu2020helper}\end{tabular}}
& F$_0$ & 0.5
\\
&
& Cr$_0$ & 0.5
\\
&
& A$_{\rm{rate}}$ & [2, 4, 6, 8]
\\
&
& H$_{\rm{m}}$ & [1, 3, 5]
\\
&
& NP$_{\rm{m}}$ & 12
\\
&
& NP$_{\rm{min}}$ & 40
\\
&
& lamda & [10, 20, 30, 40]
\\
&
& n$_0$ & [1, 2, 3]
\\
&
& gamma & [0.05, 0.1, 0.2]
\\ \cline{2-4}
& \multirow{8}{*}{\begin{tabular}[c]{@{}c@{}}MadDE\\\cite{biswas2021improving}\end{tabular}}
& p & [0.09, 0.18, 0.27, 0.36]
\\
&
& P$_{\rm{qBX}}$ & [0.01, 0.1, 0.2, 0.3, 0.5]
\\
&
& F$_0$ & 0.2
\\
&
& Cr$_0$ & 0.2
\\
&
& A$_{\rm{rate}}$ & [1.3, 1.8, 2.3, 2.8 ,3.3]
\\
&
& H$_{\rm{m}}$ & [5, 10 ,15, 20]
\\
&
& NP$_{\rm{m}}$ & [2, 4, 6, 8]
\\
&
& NP$_{\rm{min}}$ & 4
\\ \cline{2-4}
& \multirow{12}{*}{\begin{tabular}[c]{@{}c@{}}AMCDE\\\cite{ye2023differential}\end{tabular}}
& F$_0$ & 0.2
\\
&
& A$_{\rm{rate}}$ & [1.6, 2.1, 2.6, 3.1, 3.6]
\\ 
& 
& H$_{\rm{m}}$ & [5, 10, 15, 20]
\\
&
& NP$_{\rm{m}}$ & [3, 6, 9]
\\
&
& NP$_{\rm{min}}$ & 4
\\
&
& Gn & 5
\\
&
& pbc$_1$ & [0.4, 0.5, 0.6]
\\
&
& pbc$_2$ & [0.4, 0.5, 0.6]
\\
&
& pw & [0.1, 0.2, 0.3]
\\
&
& pr & [0.005, 0.01, 0.05]
\\
&
& pls$_{\rm{succ}}$ & 0.1
\\
&
& pls$_{\rm{fail}}$ & 0.0001
\\ \hline
\multirow{20}{*}{PSO}
& \multirow{3}{*}{\begin{tabular}[c]{@{}c@{}}Vanilla PSO\\\cite{pso}\end{tabular}}
& NP & [10, 20, 50, 100, 200]
\\
&
& phi$_1$ & [1, 2, 3]
\\
&
& phi$_2$ & [1, 2, 3]
\\ \cline{2-4}
& \multirow{7}{*}{\begin{tabular}[c]{@{}c@{}}GLPSO\\\cite{gong2015genetic}\end{tabular}}
& pm & [0.01, 0.1, 0.2]
\\
&
& NP & [10, 20, 50, 100, 200]
\\
&
& nsel & 10
\\
&
& w & 0.7298
\\
&
& c$_1$ & 1.49618
\\
&
& sg & 7
\\
&
& rho & [0.1, 0.2, 0.3]
\\ \cline{2-4}
& \multirow{10}{*}{\begin{tabular}[c]{@{}c@{}}sDMS-PSO\\\cite{wu2022employing}\end{tabular}}
& w & [0.729, 0.271, 0.5]
\\
&
& NP & [33, 66, 99, 198]
\\
&
& c$_1$ & [1.49445, 3., 0.75]
\\
&
& c$_2$ & [1.49445, 3., 0.75]
\\
&
& m & [1, 3, 5]
\\
&
& R & [5, 10, 15]
\\
&
& LP & [5, 10, 15]
\\
&
& LA & 8
\\
&
& L & 100
\\
&
& L\_FEs & 200
\\ \cline{1-4}
\multirow{13}{*}{PSO}
&  \multirow{13}{*}{\begin{tabular}[c]{@{}c@{}}DTPSO\\\cite{lu2023double}\end{tabular}}
& p & [0.1, 0.5, 0.9]
\\
&
& sigma & [0.25, 0.5, 0.75]
\\
&
& gamma & [0.25, 0.5, 0.75]
\\
&
& u$_1$ & [0, 0.5]
\\
&
& u$_2$ & [0, 0.5]
\\
&
& c$_{1,1}$ & [0, 1.711897]
\\
&
& c$_{1,2}$ & [0, 1.711897]
\\
&
& c$_{2,1}$ & [0, 1.711897]
\\
&
& c$_{2,2}$ & [0, 1.711897]
\\
&
& ws & 0.9
\\
&
& we & 0.4
\\
&
& NP$_{\rm{init}}$ & [50, 100, 200]
\\
&
& radius & [0.05, 0.1, 0.2]
\\ \hline
\multirow{15}{*}{ES}
& \multirow{3}{*}{\begin{tabular}[c]{@{}c@{}}SEP-CMA-ES\\\cite{ros2008simple}\end{tabular}}
& n\_individuals & [10, 20, 50, 100]
\\
&
& c\_c & [1, 2, 3, 4, 5]
\\
&
& sigma & [0.1, 0.3, 0.5]
\\ \cline{2-4}
& \multirow{6}{*}{\begin{tabular}[c]{@{}c@{}}BIPOP-CMA-ES\\\cite{hansen2009benchmarking}\end{tabular}}
& NP & [10, 20, 50, 100]
\\
&
& elite\_ratio & [0.2, 0.5, 0.7]
\\
&
& sigma\_init & 1
\\
&
& mean\_decay & 0
\\
&
& min\_num\_gens & [10, 30, 50]
\\
&
& popsize\_multiplier & [1, 2, 3, 4, 5]
\\ \cline{2-4}
& \multirow{5}{*}{\begin{tabular}[c]{@{}c@{}}MMES\\\cite{new-es}\end{tabular}}
& a\_z & [0.05, 0.1, 0.2]
\\
&
& c\_s & [0.1, 0.3, 0.5]
\\
&
& ms & [2, 4, 6]
\\
&
& n\_individuals & [25, 50, 100]
\\
&
& n\_parents & [25, 50, 100]
\\
&
& sigma & [0.1, 0.3, 0.5]
\\ \hline
\multirow{6}{*}{BO}
& \multirow{3}{*}{\begin{tabular}[c]{@{}c@{}}Vanilla BO\\\cite{bo}\end{tabular}}
& acq\_func & [LCB, EI, PI, gp\_hedge, EIps, PIps]
\\
&
& n\_initial\_points & [5, 10, 20]
\\
&
& initial\_point\_generator & [random, sobol, halton, hammersly, lhs]
\\ \cline{2-4}
& \multirow{3}{*}{\begin{tabular}[c]{@{}c@{}}LA-MCTS\\\cite{wang2020learning}\end{tabular}}
& n\_individuals & [10, 20, 50, 100]
\\
&
& c\_e & [0.01, 0.05, 0.1]
\\
&
& leaf\_size & [10, 20, 30, 40, 50]
\\ \hline
\multirow{15}{*}{LS}
& \multirow{8}{*}{\begin{tabular}[c]{@{}c@{}}Simulated Annealing\\\cite{kirkpatrick1983optimization}\end{tabular}}
& NP & [10, 20, 50, 100, 200]
\\
&
& sigma\_init & [0.1, 0.3, 0.5]
\\
&
& sigma\_decay & 1
\\
&
& sigma\_limit & [0.01, 0.05, 0.1]
\\
&
& temp\_init & 1
\\
&
& temp\_limit & 0.1
\\
&
& temp\_decay & [0.9, 0.99, 0.999]
\\
&
& boltzmann\_const & [1, 5, 10]
\\ \cline{2-4}
& \multirow{3}{*}{\begin{tabular}[c]{@{}c@{}}Dual Annealing\\\cite{xiang1997generalized}\end{tabular}}
& initial\_temp & [523, 5230, 50000]
\\
&
& visit & [1.62, 2.62, 3.62]
\\
&
& restart\_temp\_ratio & [2e-5, 2e-3, 2e-1]
\\ \cline{2-4}
& \multirow{4}{*}{\begin{tabular}[c]{@{}c@{}}NSA\\\cite{new-sa}\end{tabular}}
& sigma & [0.1, 0.3, 0.5]
\\
&
& schedule & [linear, quadratic]
\\
&
& n\_samples & [10, 20, 50, 100, 200]
\\
&
& rt & [0.9, 0.99, 0.999]
\\ \hline
\multirow{9}{*}{NO}
& \multirow{2}{*}{\begin{tabular}[c]{@{}c@{}}SLSQP\\\cite{slsqp}\end{tabular}}
& \multirow{2}{*}{eps} & \multirow{2}{*}{[1e-12, 1e-10, 1e-8, 1e-6, 1e-4]}
\\&&&
\\ \cline{2-4}
& \multirow{3}{*}{\begin{tabular}[c]{@{}c@{}}Trust-Constr\\\cite{conn2000trust}\end{tabular}}
& initial\_tr\_radius & [0.5, 1, 1.5, 2]
\\
&
& initial\_constr\_penalty & [0.5, 1, 1.5, 2]
\\
&
& factorization\_method & [equality\_constrained\_sqp, tr\_interior\_point]
\\ \cline{2-4}
& \multirow{2}{*}{\begin{tabular}[c]{@{}c@{}}COBYLA\\\cite{cobyla}\end{tabular}}
& \multirow{2}{*}{rhobeg} & \multirow{2}{*}{[0.5, 1, 1.5, 2]}
\\&&&
\\ \cline{2-4}
& \multirow{2}{*}{\begin{tabular}[c]{@{}c@{}}L-BFGS-B\\\cite{morales2011remark}\end{tabular}}
& maxcor & [5, 10, 15, 20]
\\
& 
& eps & [1e-12, 1e-10, 1e-8, 1e-6, 1e-4]
\\ \hline
\end{longtable}
\normalsize

\section{Details of Data Augmentation}\label{appx:rephrase}
\begin{figure}[h]
\vskip 0.2in
\begin{center}
\centerline{\includegraphics[width=0.9\textwidth]{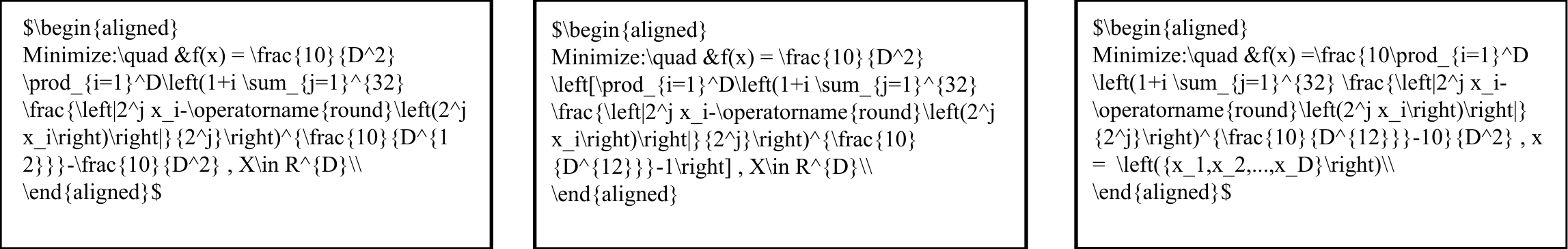}}
\caption{Three writing styles in LaTeX of the toy problem.}
\label{fig:Latex-codes}
\end{center}
\vskip -0.2in
\end{figure}

It is a common practice to augment the training data for boosting the generalization performance in recent LLMs works~\cite{formatteddata,IT2,IT-scaling}. In LLaMoCo, we alter different writing styles of a problem's definition to generate moderate diverse prompts for each problem instance generated in $P$. For the different writing styles, we investigate $500$ university students majoring in computer science, invite them to write Python or LaTeX code that they believe is correct for defining the given problem instances. After systematic statistics, we have empirically summarized several writing patterns, which we believe could approximately represent the major writing patterns of different users. Based on these different patterns, for each problem instance $p\in P$, we can obtain moderate rephrased versions for its objective function and constraints written by either Python or LaTeX code. We showcase the found patterns on a toy Katsuura problem which holds the formulation as:\\
\begin{center}
    $\begin{aligned}
Minimize:\quad &f(x) = \frac{10}{D^2} \prod_{i=1}^D\left(1+i \sum_{j=1}^{32} \frac{\left|2^j x_i-\operatorname{round}\left(2^j x_i\right)\right|}{2^j}\right)^{\frac{10}{D^{1.2}}}-\frac{10}{D^2} , X\in R^{D}\\
\end{aligned}$ \\
%
\end{center}
For LaTeX patterns, we found three different writing styles from the $500$ testees, which differ from each other mainly based on the laws of arithmetic, e.g., commutative law, distributive law and associative law. We illustrate some different LaTeX codes for our toy problem in \cref{fig:Latex-codes}.

\begin{figure}[h]
\vskip 0.2in
\begin{center}
\centerline{\includegraphics[width=0.9\textwidth]{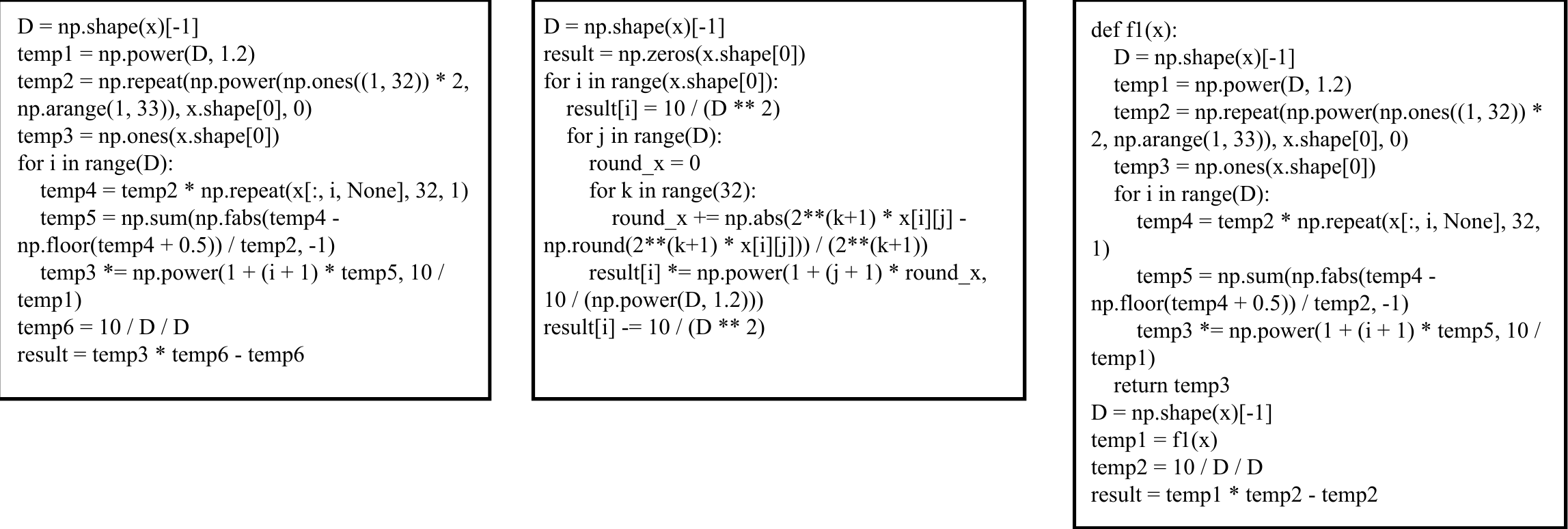}}
\caption{Three writing styles in Python of the toy problem.}
\label{fig:Python-codes}
\end{center}
\vskip -0.2in
\end{figure}

For Python patterns, the testees show different coding preferences on the writing styles of the objective functions and the constraints, e.g., some may prefer using temporary variables to store interim calculation results, some leverage \textit{numpy} to facilitate matrix operations while others use a \textit{for loop}, some may encapsulate the calculation details into a functional module etc. In \cref{fig:Python-codes} we list some of these writing styles on the toy problem.

\section{Details in Training Process}\label{appx:training_detail}
\textbf{Homogeneous batch sampling.} We further apply a homogeneous batch sampling strategy at the instruction tuning phase to reinforce the alignment of the different rephrasing version prompts for a problem $p \in P$. Concretely, we force the LLMs to sample data pairs which come from the same problem instances in a mini-batch. We observe consistent boosts in the training of LLaMoCo-S, LLaMoCo-M and LLaMoCo-L. By presenting the LLMs with a batch of homogeneous samples, they can learn patterns specific to these cross-modal prompts data more effectively.

\textbf{Batch size.} We would clarify that due to the resource limitation, all of our experiments are run on an NVIDIA A$800$ GPU. When we train the CodeGen-Mono~($350$M), the batch size is $4$ for both phases in our two-phase learning strategy. However, for one A$800$, Phi-$2$~($2.7$B) and Code Llama~($7$B) are too large to include a batch of $4$ samples, even if we adapt LoRA for them. For Phi-$2$, the batch size is $3$ and $2$ for each learning phase, while $3$ and $1$ for Code Llama.
\section{Calculation of Experimental Statistics}\label{appx:calculation}
To provide a thorough evaluation on the LLMs fine-tuned by our LLaMoCo and the other approaches, for a group of $N_p$ problem instances, we first leverage the optimization programs generated by each LLM to optimize them, for $5$ independent runs. Then we calculate the average error rate, recovery cost, optimization performance and computational cost of an approach as the performance metrics of overall performance. The calculation details of these four performance metrics in our experimental results are as follows:

\textbf{Error rate~(Err.)} The robustness of the generated optimization program is a very important index for quality of service~(QoS). We measure the robustness by the proportion of error programs generated by an LLM, named as error rate. For each instance, we use the optimization program generated by an LLM~(ours or the others) to optimize that instance for $5$ independent runs. We count the number of the generated programs which encounter compilation error or runtime error when being executed, denoted as $N_{err}$~(every single run on each instance is counted). Then the error rate of an approach on the tested instances is calculated as $\frac{N_{err}}{5\times N_p}$.

\textbf{Recovery cost~(Rec.)} While an optimization program may encounter compilation error or runtime error, we observe from our experiments that a certain proportion of the error programs could be repaired and recovered. We provide a metric named recovery cost to measure the efforts required to repair the generated programs. Concretely, for an optimization program $a_j$, we denote the number of lines in it as $L^{(j)}$, and the number of lines that need to be repaired as $L_{err}^{(j)}$. Then the recovery cost for $a_j$ is $r_j = \frac{L_{err}^{(j)}}{L^(j)}$, and the recovery cost considering all $N_{err}$ error programs is calculated as $\frac{\sum_{j=1}^{N_{err}}r_j}{N_{err}}$.

\textbf{Optimization performance~(Perf.)} We measure the optimization performance of an approach by a min-max normalized objective value descent. Concretely, we first estimate an optimal objective value $f^*_i$ for $i$-th problem instance, which can be easily obtained from our benchmarking process~(achieved best objective value). For the given approach, we denote the performance on the $i$-th problem instance in $j$-th run as a min-max normalized term $w_{i,j} = \frac{f^*_{i,j}-f^*_i}{f^0_{i,j}-f^*_i}$, where $f^0_{i,j}$ is the best objective value of the solutions initialized by the optimizer on solving the $i$-th problem instance in $j$-th run, and $f^*_{i,j}$ is the corresponding best objective the optimizer finds. Then the overall average optimization performance of the given approach on the $N_p$ instances can be calculated as follows: $\frac{\sum_{i=1}^{N_p}\sum_{j=1}^{5}w_{i,j}}{5\times N_p}$.

\textbf{Computational overhead~(Comp.)} Measuring the computational overhead by the wall-time complexity of an LLM-based approach is impractical since some of the LLMs only provide API for users. The network communication budget through calling the APIs would bias the ground results. We instead count the average number of tokens~(input+output) consumed by an approach for solving a problem instance over the test runs.
\newpage
\section{Example Input-Output of the Fine-tuned Model}\label{appx:additional}
\subsection{Synthetic unconstrained example}
\begin{figure}[h]
\centering
\includegraphics[width=0.9\textwidth]{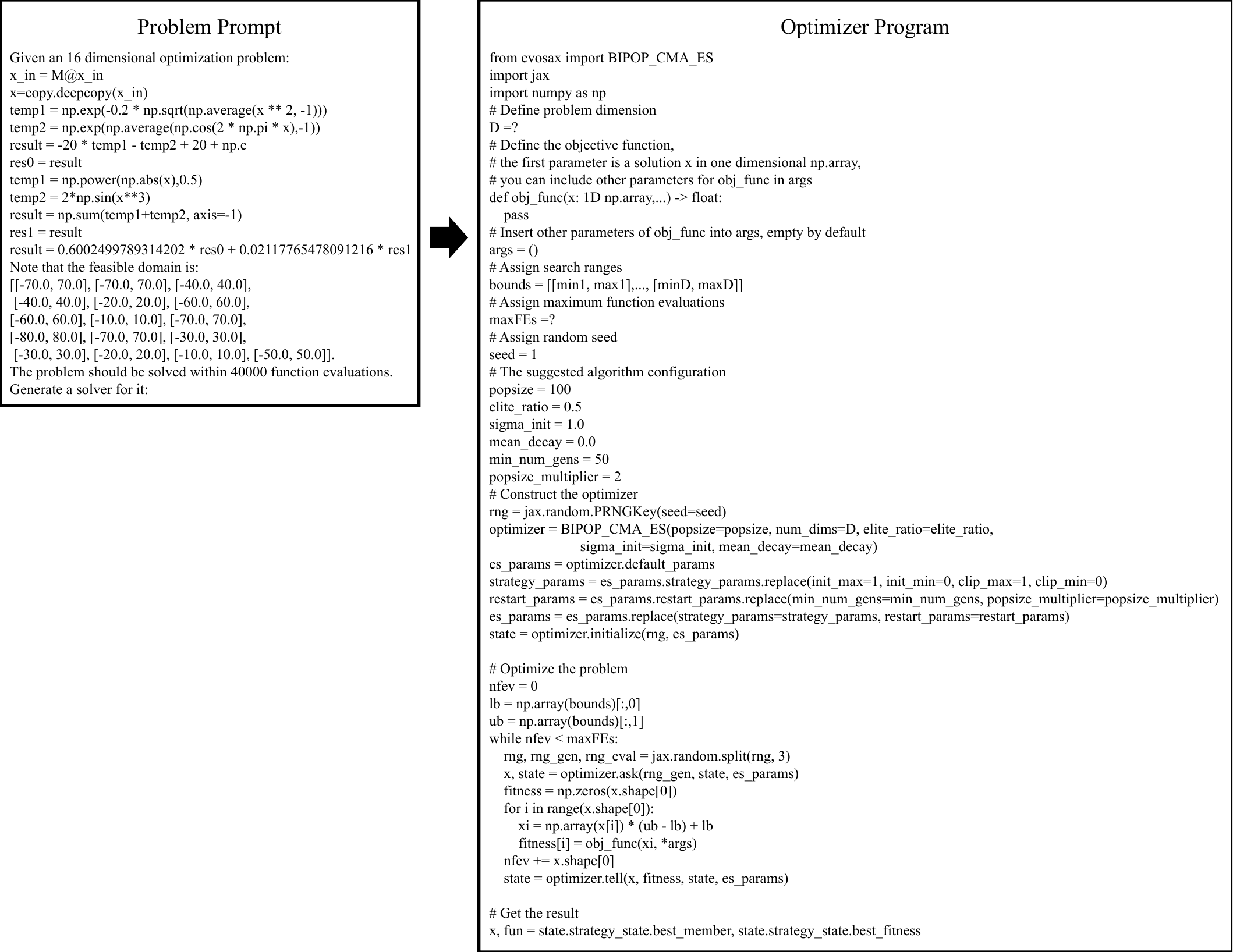}
\caption{A unconstrained problem prompt~(on the left, a Python version), and the optimization program~(on the right) output by LLaMoCo-S. The corresponding $16$-dimensional problem is constructed by a composition of two basic functions. Our LLaMoCo-S is prompted to output a competent optimizer for solving the problem within $40000$ function evaluations, which in this case, is BIPOP-CMA-ES.}
\label{fig:program1}
\end{figure}

We showcase the prompt and the generated optimization program~(\cref{fig:program1}) of a synthetic problem instance without constraints, which has the following formulation:
\begin{center}
    $\begin{aligned}
Minimize:\quad& f(x) = \sum_{i=0}^{1}W_if_{i}(z) ,z = \mathbf{M}^\mathrm{T}x , x\in R^{D} ,\mathbf{M}\in R^{D\times D}\\Where:\quad &f_{0}(\mathbf{x}) = -20 \exp \left(-0.2 \sqrt{(1/D) \sum_{i=1}^D x_i^2}\right)-\exp \left((1/D) \sum_{i=1}^D \cos \left(2 \pi x_i\right)\right)+20+e\\&f_{1}(\mathbf{x}) = \sum_{i=1}^D\left(\sqrt{\left|x_i\right|}+2sin(x_i^3)\right)\\&W0=0.6002499789314202\\&W1=0.02117765478091216\\
\end{aligned}$
\end{center}
\newpage
\subsection{Synthetic constrained example}
\begin{figure}[h]
\centering
\includegraphics[width=0.9\textwidth]{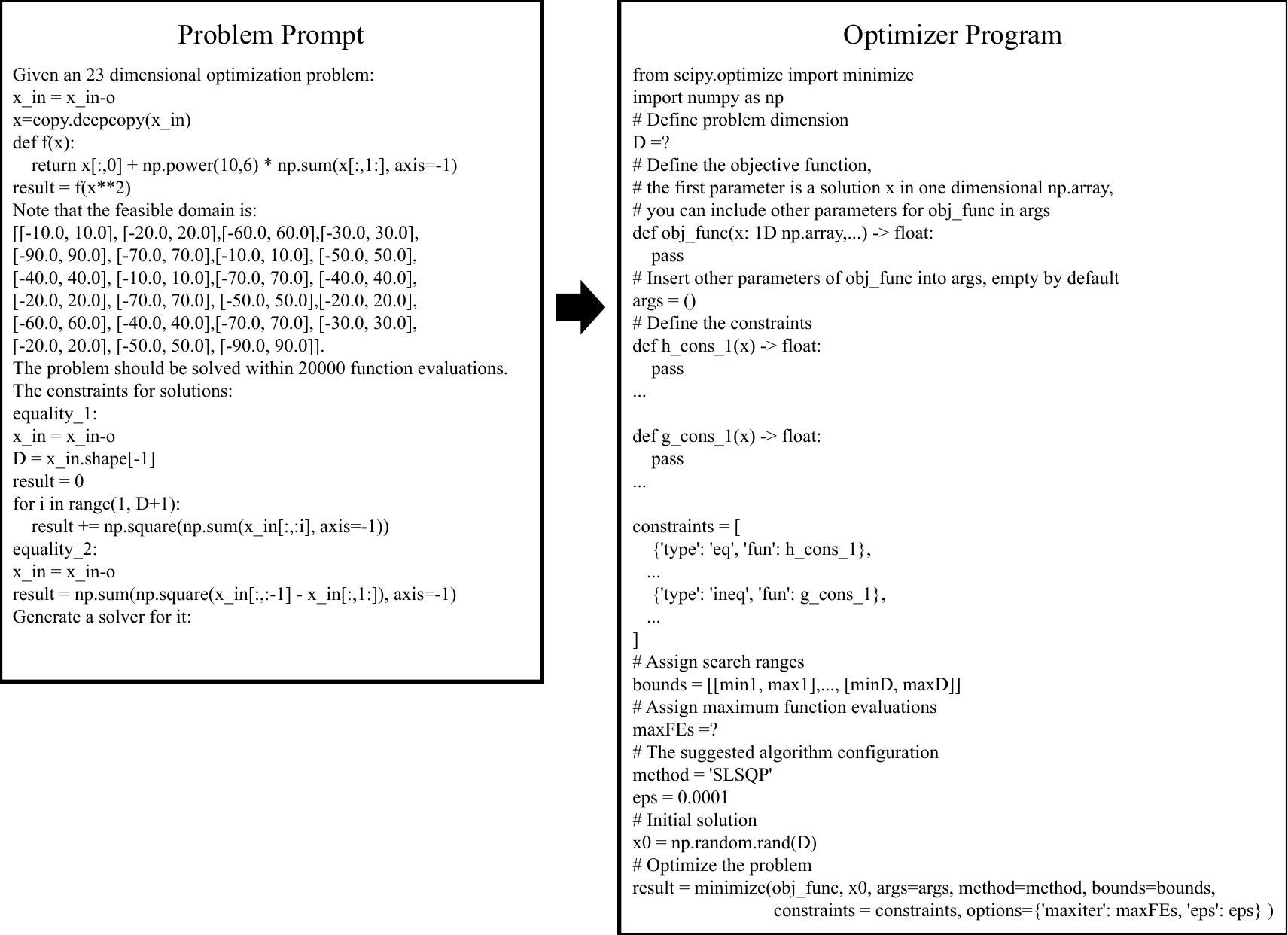}
\caption{A constrained problem prompt~(on the left, a Python version), and the optimization program~(on the right) output by LLaMoCo-S. The corresponding $23$-dimensional problem is one of the basic functions, with two additional quality constraints. Our LLaMoCo-S is prompted to output a competent optimizer for solving that problem within $20000$ function evaluations, which in this case, is SLSQP. We note that the GPT-$4$ Turbo attain the same answer on this problem. However, the configurations suggested by LLaMoCo-S achieve higher optimization performance against GPT-$4$ Turbo that adopts the default configurations.}
\label{fig:program2}
\end{figure}

We showcase the prompt and the generated optimization program~(\cref{fig:program2}) of a synthetic problem instance with some constraints, which has the following formulation:
\begin{center}
    $\begin{aligned}
Minimize:\quad & f(x) = z_1^2+10^6 \sum_{i=2}^D z_i^2 ,z = x-o , X\in R^{D} ,o\in R^{D}\\s.t.:\quad&\\ & h_0(x) : \sum_{i=1}^D\left(\sum_{j=1}^{i}y_{j}\right)^2= 0 ,y = x-o\\&h_1(x) : \sum_{i=1}^{D-1}\left(y_{i}^2-y_{i+1}\right)^2= 0 ,y = x-o\\
\end{aligned}$
\end{center}
\newpage
\subsection{Realistic example}
\begin{figure}[h]
\centering
\includegraphics[width=0.9\textwidth]{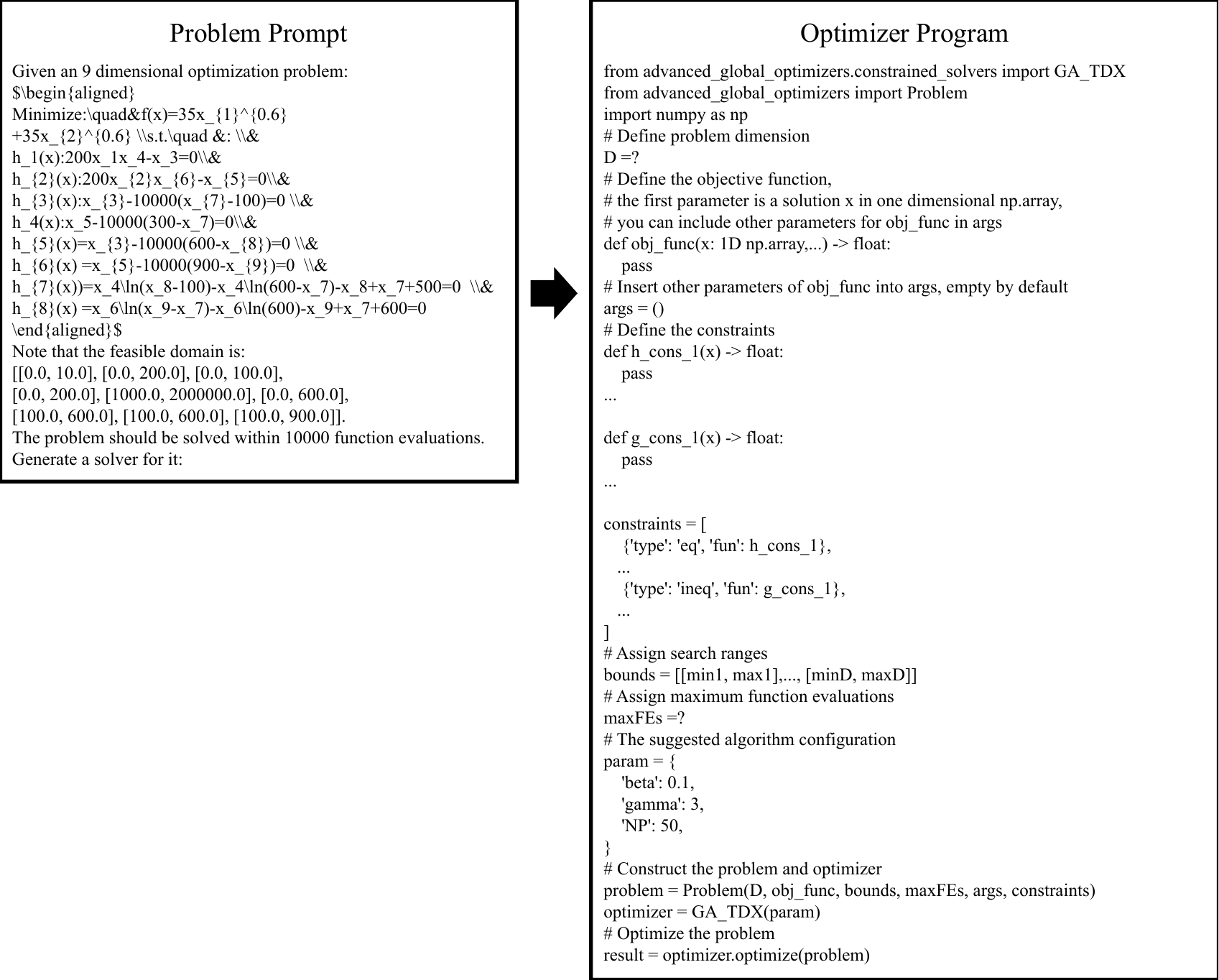}
\caption{A realistic problem prompt~(on the left, a LaTeX version), and the optimization program~(on the right) output by LLaMoCo-S. The corresponding $9$-dimensional problem holds an out-of-distribution structure, with far more constraints than the problem instances LLaMoCo-S has ever seen. Our LLaMoCo-S is prompted to output a competent optimizer for solving that problem within $10000$ function evaluations, which in this case, is an advanced GA-TDX algorithm specialized in constraints handling. We note that the GPT-$4$ Turbo suggests a DE algorithm for this problem, which is hardly adopted for solving constrianed problems.}
\label{fig:program3}
\end{figure}

We showcase the prompt and the generated optimization program~(\cref{fig:program3}) of a realistic problem instance with a large number of constraints yet with a relatively simpler objective function, which holds a different problem structure against the synthetic problems, which has the following formulation:
\begin{center}
    $\begin{aligned}
Minimize:\quad&f(x)=35x_{1}^{0.6}+35x_{2}^{0.6} \\s.t.: \quad &\\&h_1(x):200x_1x_4-x_3=0\\&h_{2}(x):200x_{2}x_{6}-x_{5}=0\\&h_{3}(x):x_{3}-10000(x_{7}-100)=0 \\&h_4(x):x_5-10000(300-x_7)=0\\&h_{5}(x)=x_{3}-10000(600-x_{8})=0 \\&h_{6}(x) =x_{5}-10000(900-x_{9})=0  \\&h_{7}(x))=x_4\ln(x_8-100)-x_4\ln(600-x_7)-x_8+x_7+500=0  \\&h_{8}(x) =x_6\ln(x_9-x_7)-x_6\ln(600)-x_9+x_7+600=0
\end{aligned}$
\end{center}

\end{document}